\documentclass[12pt,a4paper]{article}

\usepackage{pdflscape}
\usepackage{amsmath}
\usepackage{amssymb}
\usepackage{latexsym}
\usepackage{srcltx}
\usepackage{graphics}
\usepackage{color}
\usepackage{epsfig}
\usepackage{hyperref}

\newcommand{\re}{{\mathbb R}}
\newcommand{\n}{{\mathbb N}}

\newcommand{\Tn}{{\mathbb T}_n}

\newcommand{\z}{{\mathbb Z}}

\newcommand{\cF}{{\mathcal{F}}}

\newcommand{\cJ}{{\mathcal{J}}}

\newcommand{\cT}{{\mathcal{T}}}

\newcommand{\cS}{{\mathcal{S}}}

\newcommand{\cP}{{\mathcal{P}}}

\newcommand{\nill}{\mbox{\em nill}\, }

\newcommand{\bO}{{\boldsymbol{0}}}

\newcommand{\ba}{{\boldsymbol{a}}}
\newcommand{\bb}{{\boldsymbol{b}}}
\newcommand{\bc}{{\boldsymbol{c}}}
\newcommand{\bd}{{\boldsymbol{d}}}
\newcommand{\bh}{{\boldsymbol{h}}}
\newcommand{\bk}{{\boldsymbol{k}}}

\newcommand{\bp}{{\boldsymbol{p}}}
\newcommand{\bq}{{\boldsymbol{q}}}
\newcommand{\bs}{{\boldsymbol{s}}}
\newcommand{\bt}{{\boldsymbol{t}}}
\newcommand{\bv}{{\boldsymbol{v}}}
\newcommand{\bx}{{\boldsymbol{x}}}

\newcommand{\bz}{{\boldsymbol{z}}}

\newcommand{\bxi}{{\boldsymbol{\xi}}}
\newcommand{\bj}{{\boldsymbol{j}}}

\newcommand{\vardot}{\mathord{\,\cdot\,}}

\newtheorem{theorem}{Theorem}[section]
\newtheorem{prop}[theorem]{Proposition}
\newtheorem{lemma}[theorem]{Lemma}
\newtheorem{cor}[theorem]{Corollary}
\newtheorem{remark}[theorem]{Remark}
\newtheorem{ex}[theorem]{Example}

\date{}

\author{Vladimir Yu.~Protasov
\thanks{DISIM of University of L'Aquila,  {e-mail: \tt\small
vladimir.protasov@univaq.it}}
and Tatyana Zaitseva 
\thanks{Moscow State University
{e-mail: \tt\small  zaitsevatanja@gmail.com}}
}

\title{Anisotropic refinable functions and the tile B-splines 
\thanks{
The  second author is supported by by the Russian Science Foundation (project no. 23-71-30001) at Lomonosov Moscow State University.}}

\begin{document}
\maketitle

\begin{center}
{\em Dedicated to the memory of Maria Charina}  
\end{center}
\medskip 

\begin{abstract}

The regularity of refinable functions has been analysed in 
 an extensive literature and is well-understood in two cases: 1) univariate 2) 
 multivariate with an isotropic dilation matrix. The general (non-isotropic) case 
 offered a great resistance. It was done only recently by developing  the matrix method. 
 In this paper we make the next step and extend the Littlewood-Paley type method, 
 which is very efficient in the aforementioned  special cases, to general equations with arbitrary dilation matrices.  
  This gives formulas for  the higher order regularity in~$W_2^k(\re^n)$ by means of 
the Perron eigenvalue 
of a finite-dimensional linear operator on a special cone. 
Applying those results to recently introduced tile B-splines, 
we prove that they can have  a higher  smoothness than the classical ones  of the same order. 
Moreover, the two-digit tile B-splines have the minimal support of the mask 
among all refinable functions of the same order of approximation. 
This proves, in particular, the lowest  algorithmic 
  complexity of the corresponding subdivision schemes. Examples and numerical results are provided.

\smallskip


\begin{flushright}
\noindent{\bf Classification (MSCS): 65D17, 15A60, 39A99 }
\end{flushright}

\end{abstract}
\bigskip

\enlargethispage{2\baselineskip}

\section{Introduction}\label{section:intro}

\enlargethispage{2\baselineskip}

\begin{center}
\textbf{1.1. Refinement equations with matrix dilation} 
\end{center}
\smallskip 

Refinement equations are  difference functional equations with 
an integer contraction of the argument. They play a crucial role in the 
construction of wavelets in~$L_2(\re^n)$ and of the subdivision 
algorithms in the approximation theory and geometrical design, see surveys 
in~\cite{Bow, CHM1, NPS, ReifPeter08}. 
They are also applied to random power series~\cite{DDL, KM, P00n},
combinatorial number theory~\cite{FS, P17}, discrete geometry (tilings)~\cite{GroH, LW1}, etc. We consider the multivariate equation of the form:  
\begin{equation}\label{eq.ref0}
\varphi(\bx) \quad = \quad \sum_{\bk \in \z^n}\, 
c_{\bk} \, \varphi(M\bx -\bk)\, , \qquad \bx \in \re^n, 
\end{equation} 
where the dilation matrix~$M$ is integer and expanding, i.e., has integer entries and 
all eigenvalues bigger than one by modulo. We restrict ourselves to equations with 
a finite summation, i.e.,~$c_{\bk} = 0$ if~$\bk\notin Q$, where $Q$ is a finite nonempty 
subset of~$\z^n$. We also make  the standard assumption~$\sum_{\bk \in \z^n}c_{\bk} \, =\, m, 
$ where $m=|\det M|$.
In this case, equation~(\ref{eq.ref0}) possesses a unique, up to multiplication by a constant, 
 solution~$\varphi$ ({\em refinable function}) in the space~$\cS_0'(\re^n)$ of compactly-supported tempered distributions. 
In most of applications this solution needs to be at least from~$L_2(\re^n)$.  
Often a higher smoothness is required. 
Let us note that~$\varphi$ is never infinitely smooth (in~\cite{WX} this was proved for 
more general equations). 

 In the univariate case ($n=1$) the 
smoothness of refinable functions
has been studied in great detail. Two methods 
are usually most prominent. The matrix approach~\cite{CDM, CHM2, D, MR, P06} finds the biggest~$k$
such that~$\varphi\in C^k(\re)$ and computes the precise value of the 
H\"older exponent in this space. This  method  requires computation of the joint spectral radius of special matrices~$\{T_i\}_{i=0}^{m-1}$, which is a hard problem even in low 
dimensions. 
Recent tools elaborated in~\cite{GP1, GP2, M} 
 can do this for matrices of non very large sizes, which makes it possible to 
compute the regularity for univariate 
equation~$\varphi(x)  =  \sum_{k =0}^N\, 
c_{k} \varphi(2x -k)$ of order $N \le 25$. The matrix method is also applicable to 
compute the regularity in~$L_2(\re)$~\cite{LMW} and in the Sobolev spaces~$W_2^k(\re)$
(of absolutely continuous functions~$f$ such that~$f^{(k)} \in L_2$). In this case the joint spectral 
radius computation can be replaced by computing the usual spectral radius of one matrix~$T$,  
although, of quite a big size, around~$\frac12\, N^2$, see~\cite{P97}. 

The second approach is usually referred to as the {\em Sobolev regularity computation}
or the {\em Littlewood-Paley type method}. It finds  the regularity in~$W_2^k(\re)$
by means of computing  the  Sobolev regularity defined as 
$\sup\, \bigl\{ \beta \ge 0 \, : \,  \int_{\re}|\varphi(\xi)|^2 (1 + |\xi|^2)^{\beta}\, dt\, < \, \infty\bigr\}$. This problem is reduced to finding the usual spectral radius of a matrix of a much lower size, around~$2N$~\cite{CD, DD, E, V}. The shortcoming of this 
method is that it deals with  the regularity not in~$C$ but in~$L_2$. 
Say,   to prove the continuity 
of~$\varphi$, one needs to show that~$\varphi \in W^1_2$. The matrix method and the 
Littlewood-Paley type method are both 
 widely used. 

For multivariate refinement equations, the situation is more complicated. The smoothness of solutions has been computed 
only in the special case when the matrix~$M$ is {\em isotropic}, i.e., 
is similar to a diagonal matrix with equal by modulus eigenvalues. Isotropic matrices 
expand all distances equally in all directions, i.e., 
the sequence~$\{M^k\bx\}_{k\in \n}$ has the same asymptotic growth 
for all~$\bx \in \re^n\setminus\{0\}$. In this case many results on the univariate equations 
can be generalized to~$\re^n$. In particular,  both of the aforementioned methods work for isotropic dilations, although the 
dimensions of the corresponding matrices~$\{T_i\}_{i=1}^{m-1}$ and~$T$ become much higher. 
This was done in~\cite{CHM1, CDM, CJR, CD, HJ, H1, H3, J, JZ, KPS, RonS}.  
If $M$ is not isotropic, then none of the methods is applicable. 
Already simple two-dimensional examples show that the formulas 
for the H\"older exponent are not valid in the non-isotropic case. 
Even some basic questions, for instance, whether or not~$\varphi$ 
possesses a derivative (say, from~$C(\re^n)$ or from~$L_2(\re^n)$),   
remain unanswered. Here it should be emphasized that 
refinement equations with non-isotropic dilations are important in 
many problems of signal processing, wavelets and geometric design, see~\cite{ACJRV, Bow, CDR1, CDR2, CM, CP, KPS}. In~\cite{Bow, C, CGV, SF} the authors compute regularity in special anisotropic H\"older and Sobolev spaces.

It was not before 2019 that the modification of the matrix method
was elaborated for non-isotropic dilations~\cite{CP}. This gave 
closed formulas for the H\"older exponent and for higher order regularity  
of general refinable functions by means of 
joint spectral radii. However, because of the large size of the matrices~$T_i$, this 
method is applicable  mostly for very low dimensions~$n$ 
and for small sets~$Q$ of nonzero coefficients~$c_{\bk}$. 

An improvement would be achieved if one could make the next step: to modify the 
Littlewood-Paley type method for non-isotropic matrices.  In this paper we make this step.  
 We obtain a sharp criterion which, for an arbitrary dilation matrix~$M$,  
ensures that the solution of refinement equation~(\ref{eq.ref0}) belongs to~$W^k_2(\re^n)$ 
and presents a formula for its H\"older exponent. This formula is different and more 
complicated than one for  the isotropic case, although, for many equations, it can be simplified.  Which is most important, it allows us to solve the smoothness problem in higher dimensions~$n$ and 
for larger sets~$Q$. Numerical results in Section~\ref{section:numeric}
demonstrate that indeed the Littlewood-Paley method for computing the $L_2$-regularity 
of multivariate refinement equations is more efficient than the matrix approach 
because of a much smaller dimension of the problem of the spectral radius computation.

Then we apply this formula to 
the tile B-splines and to the corresponding subdivision schemes.

 \bigskip 

\begin{center}
\textbf{1.2. The tile B-splines}  
\end{center}
\bigskip

In the second part of the paper we analyse the regularity of the tile B-splines
recently introduced in~\cite{Z1}. They are analogues of the 
univariate cardinal B-splines. The cardinal B-spline of order  $\ell \ge 0$
is $B_{\ell} \, = \, \underbrace{\chi_{[0,1]} * \cdots * \chi_{[0,1]}}_{\ell + 1}$, 
it has  integer nodes, is supported on~$[0, \ell+1]$, and 
belongs to~$W_2^\ell(\re^n)$. 
In particular~$B_0 = \chi_{[0,1]}$.  Cardinal B-splines have countless applications due to  
their simplicity and many optimal properties (the smoothest refinable function of a given order, etc.)  The multivariate B-splines are defined in the standard way as 
 direct products of univariate ones. 
This approach, however, has significant shortcomings, especially in the construction of wavelets and Haar systems, see~\cite{CHM2, C, KPS} for general discussion of this aspect. In contrast, the tile B-splines 
are not obtained by direct products but by the convolutions: 
$\varphi_{\ell} \, = \, \chi_{G} * \cdots * \chi_{G}$, where 
$G$ is a self-similar tile in~$\re^n$. Here we need to recall some  definitions. 

Given an integer expanding $n\times n$ matrix~$M$, we take a
system of ``digits'' $D=\{\bd_0, \ldots , \bd_{m-1}\}$, which are $m = |\det M|$ arbitrary 
representatives of the  quotient group~$\z^n /M\z^n$. The set~$G$ of $M$-adic 
fractions: 
$$
G\ = \ \Bigl\{0.\bd_{i_1}\bd_{i_2}\cdots \  = \ 
\sum_{k\in \n}\, M^{-k}\bd_{i_k} \ : \  \bd_{i_k} \in D\, , \ k\in \n
 \Bigr\}
$$
is a {\em self-similar attractor}. This  is a compact set of an integer positive Lebesgue measure~\cite{GroH, LW1}. 
If its measure is equal to one, then the attractor is called {\em a tile}. 
Integer shifts of a tile form a partition of~$\re^n$. 
Moreover, a tile is self-similar: it is a union of several parallel shifts of the 
set~$M^{-1}G$. Its indicator~$\varphi_0 = \chi_{G}$
satisfies the refinement equation~$\varphi_0(\bx)\, = \, \sum_{i=0}^{m-1}\varphi_0(M\bx - \bd_i) $. 
In some sense, the tile is a unit segment in the~$M$-adic system on~$\re^n$. 
The Haar systems in~$L_2(\re^n)$ with the dilation factor~$M$ can be obtained from the function~$\chi_{G}$
by the same procedure as the classical Haar system with the double dilation in~$L_2(\re)$
is obtained from~$\chi_{[0,1]}$. 
To construct wavelet bases in~$L_2(\re^n)$ 
one can naturally define the tile B-spline~$\varphi_{\ell}$ as 
a convolution of $\ell+1$ functions~$\varphi_0$. Strictly speaking, $\varphi_{\ell}$
is not a spline, since it is not a piecewise-polynomial function. 
On the other hand, $\varphi_{\ell}$ is smooth, generates algebraic 
polynomials of degrees~$\ \le \ell$ by its integer translates, and the  values~$\varphi_{\ell}(\bx)$ 
 are  easily computed  in all $M$-adic points~$\bx$
(from  the corresponding refinement equation). The tile~B-splines are  useful in approximations and interpolations 
as well as for constructing of multivariate wavelets and subdivision schemes~\cite{Z1}. 

Applying our main results, we prove that  some  families of tile~B-splines have  
(quite surprisingly) a higher  smoothness than the standard B-splines  of the same order. 
Note that nothing of this kind occurs in the univariate case, where 
 the cardinal B-spline is the smoothest refinable function of a given order. 
 This  ``super-smoothness'' yields, in particular, that 
 the subdivision schemes defined by tile B-splines generate smoother surfaces and possess faster convergence. 
 Moreover, in Section~\ref{section:optimality} we prove the optimality of tile B-splines among 
 all refinable functions of a given order.

 \bigskip 

\begin{center}
\textbf{1.3. The tile subdivision schemes}  
\end{center}
\bigskip 

The tile B-splines demonstrate their efficiency in the construction 
of multivariate wavelets~\cite{Z1} and of subdivision schemes. 
The {\em subdivision operator}  $S: \,  \ell_{\infty}(\z^n) \to 
\ell_{\infty}(\z^n)$
generated  by a dilation matrix~$M$ and coefficients~$\{c_{\bk}\}_{\bk \in \z^n}$ 
is defined as~$\, [S\ba](\cdot) \, = \, \sum_{\bk \in \z^n} c_{\cdot - M\bk}\, \ba(\bk)$. 
The subdivision scheme is a recursive  application of~$S$ to  a function~$f_0: \z^n \, \to \re$. The subdivision scheme {\em converges 
to a limit function}~$f \in C(\re^n)$ if 
$\bigl\|f(M^{-j} \cdot )  \, - \, S^jf_0(\cdot )\bigr\|_{\infty} \to \infty$ as~$j\to \infty$. 
Similarly one defines the convergence in~$L_2$. In Section~\ref{section:subdivisions}
we analyse the properties of  subdivision schemes generated by the tile B-splines, in particular, those having the super-smoothness.  
We prove that the tile schemes are optimal in the following sense: 
they have the lowest possible complexity of one iteration 
among all schemes of a given approximation rate.

 \bigskip 

\begin{center}
\textbf{1.4. Notation}  
\end{center}
\bigskip 

We denote vectors by bold letters and scalars by usual ones, 
so~$\bx \, = \, ( x_1, \ldots , x_n) \in\re^n$. The polynomials are also denoted 
by bold letters. The cardinality of a finite set~$X$ is denoted by~$|X|$. 

We define the Fourier transform in~$\re^n$ by the formula
 $\widehat{f} (\bxi)\, =\, \int_{\re^n}f(\bx)e^{-2\pi i (\bx, \bxi)}d\bx$. 
 As usual, the convolution of functions $f, g$ is defined by
$f * g (x) \ = \ \int f(t) g(x - t) dt$, also we  
use the  symmetrized autoconvolution 
of~a function~$\varphi$ defined as  $\varphi * (\varphi\_)$, where $\varphi\_(\vardot) = \varphi (-\vardot)$. 
 
We consider trigonometric polynomials of the 
form~$\bp(\bxi)\, = \, \sum_{\bk \in \Omega} c_{\bk}e^{-2\pi i (\bk, \xi)},\, \bxi \in \re^n$, where $\Omega$ is a finite subset of~$\z^n$. We  usually drop the word ``trigonometric''. 
The set of zeros of a function~$f$ will be denoted as  $\nill (f)\, = \, \{\bx \in \re^n : \ f(\bx) = 0\}$.  Similarly, for a natural number~$\ell$, $\nill_{\ell}(f)$ denotes the set of zeros of order at least~$\ell+1$. 
For any closed subset $C \subset \re^n$, let  $\pi\{C\}$ be the set of points 
$\bxi \in \re^n$ 
such that 
$\bxi + \bk \in C$ for all~$\bk \in \z^n$. 
For several subsets~$C_i$, we 
write~$\pi\, \{C_1, \ldots , C_m\}$ for 
$\pi\{C_1 \cup  \ldots \cup C_m\}$. 

Two functions $f(x), g(x)$ on~$\re$, 
are asymptotically equivalent $\, \bigl(f(x)  \ {\asymp} \ g(x)\bigr)$ if 
there exist positive $c_1, c_2$ such that 
$\, c_1|f(x)|\, \le \, |g(x)| \, \le \, c_2|f(x)|, \ x \in \re$.
We also define the {\em exponential  asymptotic equivalence}. 
For a given sequence  $\{a_k\}_{k=1}^{\infty}$ and $\lambda > 0$, 
we write $a_k  \ \dot{\asymp} \ \lambda^k$  if 
 there exist positive $r, c_1, c_2$  such that 
$$
 c_1\lambda^k\ \le \ |a_k| \ \le \ c_2 k^r \lambda^k, \ k \in \n.
 $$ 
We write $|a_k|  \ \dot{\le} \ \lambda^k$ if there exist positive constants $r, c$
such that~$|a_k|  \ \le \ c \, k^r \, \lambda^k, \, k\in \n$.  

For a given subset~$X \subset \re^n$, we denote by ${\rm cone}\, X \, = \, {\rm co}\, \{\lambda \bx \ : \ \bx \in X, \, \lambda \in \re_+\}$ 
 the conic hull of~$X$. 

The function $\bc(\bxi)\, = \, \sum_{\bk \in \z^n} c_{\bk}e^{-2\pi i (\bk, \bxi)}$ is the {\em mask}
of the refinement  equation~(\ref{eq.ref0}). We say that the mask 
{\em is supported on a set~$\Omega \subset \z^n$} if $c_{\bk} = 0$ for~$\bk \notin \Omega$. 
Since we assume that the support is finite,  $\bc(\bxi)$ is a trigonometric polynomial.

\bigskip

\bigskip 

\section{The roadmap for the main results}\label{section:road}
\bigskip 

\begin{center}
\textbf{2.1. Computing the regularity in the isotropic case: an overview}
\end{center}
\bigskip 

The matrix approach makes use of the  {\em transition operator}: 
$$
[\cF f](\vardot)\quad = \quad \sum_{\bk \in \z^n} c_{\bk} f(M\vardot \, - \, \bk),
$$
 where, recall, 
the summation is finite and the mask is supported on the subset~$\Omega\subset \z^n$. 
The refinable function~$\varphi$ is the eigenfunction of~$\cF$ with the eigenvalue one. 
Assume that $M$ is isotropic (including the case~$n=1$), 
$r=\rho(M)$ is the spectral radius of~$M$,~$W$ is the linear space spanned by the differences  $\varphi(\vardot +\bh) - 
\varphi (\vardot), \, \bh \in \re^n$. From the refinement equation it 
follows that the value~$\sup_{\|\bh\| \le r^{-k}}\|\varphi(\vardot +\bh) - 
\varphi (\vardot)\|_{\infty}$ is equivalent to~$\|\cF^k|_{W}\|$. 
This implies that the H\"older exponent of $\varphi$
is equal to~$\log_{\,1/r} \, \rho(\cF|_{W})$~\cite{CDM, CHM1}. Let us stress that this is true only if~$M$ is isotropic.  
It is proved that the spectral radius $\rho(\cF|_{W})$ in the space~$C(\re^n)$ is equal to 
the so-called joint spectral radius of $N\times N$ {\em transition matrices}~$\{T_0, \ldots , T_{m-1}\}$
constructed by the coefficients~$c_{\bk}$ and restricted to a certain subspace of~$\re^N$. 
The computation of the joint spectral radius, however,  is a hard problem.   The H\"older exponent in~$L_2$
is expressed by the {\em $L_2$-radius} of the same matrices. 
It is equal to the square root of the usual spectral radius of the matrix 
$T = \frac1m \, \sum_{i=0}^{m-1} \, T_i \otimes T_i$, where~$\otimes$ denotes the Kronecker product~\cite{P97}. Thus, the computation of the~$L_2$-regularity 
by the matrix method involves the usual spectral radius~$\rho(T)$. 
However, $T$ has the size of order~$\frac12\, N^2$ while the matrices~$T_i$
have size of order~$N$. 

We consider the Littlewood-Paley type method (computation of Sobolev regularity). 
Its idea is to associate to every function~$f\in W$ the periodization of the square of its Fourier transform~$\widehat f$: 
\begin{equation}\label{eq.Phi0}
\Phi_f\, (\bxi)\quad = \quad  \sum_{\bk \in \z^n} \bigl| \widehat{f}(\bxi + \bk) \bigr|^2\, . 
\end{equation} 
The Plancherel equality yields that $\|\Phi_f\|_1 \, = \, \|\widehat f\|_2^2 \, = \, 
\|f\|_2^2$ (the norm of $\Phi_f$ is computed in $L_1(\Tn)$). If we apply the operator~$\cF$ to the function~$f$, 
what happens then to~$\Phi_f$\, ? To answer this question, we define the operator~$\cT$
on~$L_2(\Tn)$ as follows. Let $D^*\, = \, \{\bd_i^*\}_{i=0}^{m-1}$ be a set of digits 
for the transpose matrix~$M^T$.  This set may not coincide with~$D$.  Then 
\begin{equation}\label{eq.T0}
\bigl[\cT \bp \bigr](\bxi)\quad = \quad 
\sum_{j=0}^{m-1}\, \bigl| \, \bc \, (M^{-T}(\bxi + \bd^*_j)) \bigr|^2 \, \bp \bigl(\, 
M^{-T}(\bxi + \bd^*_j)\, \bigr)\ , \qquad \bp \in L_2(\Tn)\, . 
\end{equation} 
A direct calculation shows that~$\, \Phi_{\, \cF f} \, = \, \cT\, [\Phi_f]$. 
Hence, for each~$k\ge 1$, the norm~$\|\cF^k f\|_2$, is equal to the square root of the $L_1$-norm of 
$\cT^k[\Phi_f]$. One may wonder what is the use of changing one  
infinite-dimensional operator~$\cF$ to another operator~$\cT$? 
 Observe that, for every compactly supported~$f\in W$, 
the function~$\Phi_{f}$
is a polynomial. Indeed, $\Phi_{f}$ has the following Fourier expansion: 
 $$
 \Phi_{f}(\bxi) \ = \ \sum_{\bk \in \z^n}\, 
\Bigl(\, f(\vardot)\, , \, f(\vardot + \bk) \, \Bigr)\, e^{\, 2\pi i \, (\bk, \bxi)}\, , 
$$ 
where  only finitely many coefficients~$\bigl(f(\cdot), f(\cdot + \bk) \bigr)$
are nonzero. This argument reduces, the operator~$\cT$
to a finite-dimensional space of polynomials~$\cP_{\Omega}$ with coefficients supported
 on certain finite symmetric subset~$\Omega \subset \z^n$. It is possible to prove that 
 if~$f\in W$, then the polynomial $\Phi_f$ has a root of multiplicity two at zero. 
 Denoting the space of such polynomials by~$\cP_0 = \cP_0(\Omega)$, we have  
 $\rho^2(\cF|_{W})\, = \, \rho(\cT|_{\cP_0})$, and therefore, 
 the  H\"older exponent of $\varphi$ in~$L_2$
is equal to~$ \log_{1/r} \, \rho(\cF|_{W})\, = \,  \frac12\, 
 \log_{1/r} \rho(\cT|_{\cP_0})$. Note that the dimension of the space~$\cP_{\Omega}$ is 
 usually much less than  the size of the matrix~$T$, which  makes this 
 formula more efficient in practice and applicable for a wider set of refinement equations. 
 This holds, however, only for the isotropic case. The situation 
 in the case of general 
 matrix~$M$ is quite different.

\bigskip 

\begin{center}
\textbf{2.2. What is wrong with the non-isotropic case?}
\end{center}
\bigskip

The advantage of the Littlewood-Paley type method is that it expresses the exponent of 
Sobolev and of $L_2$ 
regularity merely by the spectral radius of one matrix, whose size is usually not too large.  Neither the joint spectral radius
nor the Kronecker lifting is involved.  However, the 
formula $\, \frac{1}{2}\log_{1/r}\, \rho(\cT|_{\cP_0})$ is valid only 
for isotropic~$M$. For general $M$, neither this formula, nor 
its possible 
modifications  (for example, by replacing~$r$ with $|{\rm \det}\, M|^{1/n}$) holds. 
The reason is that a non-isotropic matrix~$M$ expands the distances differently in different  
directions. This aspect has to be taken  into account in  any  regularity analysis of 
non-isotropic refinement equations. Can the Littlewood-Paley type method be 
applied to arbitrary dilation matrices? 
\medskip

\smallskip

\bigskip 

\begin{center}
\textbf{2.3. Our solution for the non-isotropic case: the four steps}
\end{center}
\bigskip

 \textbf{The  idea} of handling the non-isotropic case was suggested in~\cite{CP}
and worked out in the space~$C(\re^n)$. One considers the {\em spectral subspaces}
$J_1, \ldots , J_q$ of the matrix~$M$. Each subspace~$J_s$ corresponds 
to all eigenvalues of~$M$ of a given modulus~$r_s$. Thus, 
$r_1 > \cdots > r_q$ are all moduli of eigenvalues of~$M$ and 
the restriction $M|_{J_s}$ has all the eigenvalues of modulus~$r_s$. 
Note that  the operator~$M|_{J_s}$ is not necessarily isotropic since it may have 
nontrivial Jordan blocks. We have~$\re^n \, = \, J_1 \oplus \cdots \oplus J_q$
(the direct sum). 
The regularity of functions on~$\re^n$ is considered separately 
on the directions of the subspaces~$J_s$. 
This is done by computing the spectral radii 
of the restrictions of the operator~$\cF$ to the subspaces~$W_s$ 
spanned by the differences  $\varphi(\vardot +\bh) - 
\varphi (\vardot), \, \bh \in J_s$. 
 Then the regularity on the whole~$\re^n$ is 
 the minimal regularity along~$J_s, \, s = 1, \ldots , q$. 
 This idea works  for the matrix method~\cite{CP}, however, its extension to the Littlewood-Paley type method meets considerable difficulties.  
One of the reasons is that the map: $\Phi : \, W \to \cP_0$, 
which associates~$\Phi_f$ to a function~$f\in W$, is nonlinear.
Hence, the images of the subspaces~$W_s$ become nonlinear manifolds in~$\cP_0$.
We handle this problem by a new construction realized in four  steps described below. 
\medskip 

 \textbf{Step 1.} We use the spectral properties of linear operators 
with an invariant cone.  
Theorem~\ref{th.45}, the main result of Section~\ref{section:cone}, states that if $K\subset \re^n$ is a closed solid pointed cone respected by a linear operator~$A$, then 
for every point~$\bx \in K$, the asymptotic growth of the trajectory~$\{A^j\bx\}_{j=1}^{\infty}$ 
as~$j\to \infty$ is determined by the spectral radius of~$A$ restricted to the 
minimal face of~$K$ invariant with respect to~$A$ and containing~$\bx$. 
The same is true not only for one point~$\bx$ but also for an arbitrary subset~$\Gamma \subset K$. 
The fastest asymptotic growth of trajectories~$\{A^j\bx\}_{j=1}^{\infty}$
over all points~$\bx \in \Gamma$ is defined by the minimal invariant face of~$K$
containing~$\Gamma$.

 \medskip 

 \textbf{Step 2.} Theorem~\ref{th.45} is applied for the
 cone~$\cP_+ = \cP_+(\Omega)$
of nonnegative trigonometric polynomials $\bp(\bx)\, = \, \sum_{\bk \in \Omega}\, p_{\bk} 
\, e^{2\pi i \, (\bk, \bx)}$. To this end, we characterize all invariant faces of the 
cone~$\cP_+$. Propositions~\ref{p.60} and~\ref{p.70} assert that, for an arbitrary operator~$A$ leaving~$\cP_+$ invariant, every invariant face~$L\subset \cP_+$  is defined 
by a set~$S\subset \re^n$ of common zeros, i.e., $\ L = \bigl\{\bp \in \cP_+ \ : \ \bp(S) = 0\bigr\}$.  
\medskip 

 \textbf{Step 3.} We 
 consider an arbitrary  difference~$f=\varphi(\vardot + \bh) - \varphi(\vardot)$ along 
 a shift vector~$\bh\in J_s$. 
 Applying contractions by the powers of~$M^{-1}$, 
we obtain a small vector~$M^{-j}\bh$ of the  norm 
equivalent to~$r_s^{-j}$. On the other hand, the corresponding 
function~$\Phi_f$ after changing the shift vector  from~$\bh$
to~$M^{-j}\bh$, becomes~$\cT^j\Phi_f$, where the linear operator~$\cT$
is given by~(\ref{eq.T0}). Now we apply the results of Steps 1, 2 to the 
positive polynomial cone~$\cP_+$ and conclude that
$\|\cT^j\Phi_f\| \ \dot{\asymp}\  (\rho_s)^j$, where 
$\rho_s$ is the spectral radius of~$\cT$ restricted to a certain face~$K_s$ 
of~$\cP_+$. Thus, $\| \varphi(\vardot + M^{-j}\bh) - \varphi(\vardot)\|^2 \, = \, 
\ \|\cT^j\Phi_f\|_1 \ \dot{\asymp} \ (\rho_s)^j$. 
We see that while the norm of the shift $\bh \in J_s$ is multiplied by (roughly)~$r_s^{-j}$, 
the norm of the variation of~$\varphi$ is multiplied by~$(\rho_s)^j$. 
Now it is  easily deduced that the H\"older exponent of~$\varphi$ along the 
subspace~$J_s$ is equal to~$\, \log_{1/r_s} \rho_s$.
The minimum of those values over~$s=1, \ldots , q$ gives the  H\"older exponent of~$\varphi$. 

\smallskip

Note that if the dilation matrix~$M$ is isotropic, then there is only one spectral subspace~$J_1=\re^n$
and hence~$W_1 = W$.  Therefore,~$\Phi$ maps $W_1$ to the entire space~$\cP_0$, 
and the H\"older exponent of~$\varphi$ is expressed 
merely by the spectral radius of~$\cT|_{\cP_0}$. 
This was the core of the method in the isotropic case~\cite{CD,  E, H3, JZ}.  
Now, dealing with an arbitrary matrix~$M$, 
 we may have proper subspaces~$W_1, \ldots , W_q$
of~$\re^n$. 
The image of~$W_s$ by the map~$\Phi$ is some nonlinear
 manifold in~$\cP_0$.  If one follows the line of reasoning from the
 isotropic case, he needs to apply the operators~$\cT^j$ to those 
 manifolds, which can be difficult. We resolve this problem  by means of 
 invariant cones. Steps~1 and~2 allow us to 
 to apply the operator~$\cT^j$ to the minimal invariant face~$K_s$ of the 
cone~$\cP_+$ containing those manifolds. This minimal face~$K_s$  can be found explicitly, 
which leads to the formula for the H\"older exponent in~$L_2$ (Theorem~\ref{th.50}).

\medskip 

 \textbf{Step 4.} To extend our approach to the higher order regularity 
we use  the $k$th order differences
of the function~$\varphi$ along the shift  vectors~$\bh_i \in \re^n$. 
Replacing the  difference~$\varphi(\vardot + \bh) - \varphi(\vardot)$ 
by the higher order difference, we adopt 
our approach to the higher order regularity and obtain an explicit formula for 
the H\"older exponents of the derivatives in Theorem~\ref{th.50h}. 
 \medskip

 \bigskip 

\begin{center}
\textbf{2.4. Applications to tiles, B-splines, and subdivisions}
\end{center}
\bigskip

 We use our method (Theorems~\ref{th.50} and~\ref{th.50h})
to find the regularity of the tile B-splines. To this end, one needs to find 
the invariant faces~$K_s$ of the cone~$\cP_+$. This is done 
in Section~\ref{section:splines} by involving the geometry of cylinders in~$\re^n$ composed by tiles and some known facts on the rational subspaces of~$\re^n$ (Theorem~\ref{th.60}). 
We shall see that the~faces~$K_s$ in case of tile B-splines  have a simple structure, which leads to the precise values of the high order H\"older exponents
(Theorem~\ref{th.65}). In particular, 
we find several types of the tile B-splines that have the regularity higher 
than the traditional cardinal B-splines and the faster convergence of the corresponding subdivision schemes.  The efficiency of the tile subdivisions 
is demonstrated in Section~\ref{section:subdivisions} along with other 
examples and numerical results.

\section{Auxiliary results. Operators leaving a cone invariant} \label{section:cone}

\smallskip 

We consider a convex cone~$K \subset \re^d$, which is a set 
closed with respect to addition and multiplication by positive numbers: 
if $\ba, \bb \in K$, then $\, \ba+ \bb \in K$ and $\lambda \ba \in K$
for every~$\lambda \in \re_+$. In what follows we 
 always assume that a cone is closed, {\em solid} (or {\em full-dimensional}, i.e., 
 possesses a nonempty interior) and 
 {\em pointed}   (does not contain a straight line). 
 A linear  operator~$A$ in~$\re^d$ {\em leaves~$K$ invariant} if 
 $\ AK \subset K$.  By the Krein-Rutman theorem~\cite{KR},
an operator leaving a cone~$K$ invariant  possesses a leading eigenvector 
$\bv \in K, \, \bv \ne 0$, for which~$A\bv = \lambda_{\max}\bv, \ \lambda_{\max} = \rho(A)$. 
This fact extends the well-known  Perron-Frobenius theorem from the positive orthant~$\re^d_+$ 
to an arbitrary cone~$K \subset \re^d$.  

A plane of support of a cone~$K$ is a hyperplane that intersects~$K$  but 
does not intersect its interior. A {\em face} of a cone is 
its intersection with a plane of support. In particular, the apex~$\bO$ is the face of dimension zero. In addition, the cone~$K$ itself is considered as the face  of dimension~$d$.  
All other faces are
  {\em proper} and have dimensions from~one (the  edges) to~$d-1$
(the facets). All proper faces are cones of smaller dimensions.  A {\em face plane} is a linear span of a face.  For  an arbitrary subset~$\Gamma \subset K$, its 
{\em minimal face} is the smallest by inclusion face of~$K$ that contains~$\Gamma$. 
 If an operator~$A$ leaves~$K$ invariant, then the {\em minimal invariant face for~$\Gamma$} is the smallest by inclusion face of~$K$ invariant with respect to~$A$
and  containing $\Gamma$.

We are interested in the asymprotic behaviour of the trajectory~$A^j\bx$ as $j\to \infty$,  
for a given point~$\bx$. Clearly, $\|A^j\bx\| \, \dot{\le} \, \rho^{\, j}$, where 
$\rho = \rho(A)$ is the spectral radius of~$A$. If $\bx$ belongs to an invariant cone, one can say more:  
\begin{prop}\label{p.20}
Let an operator~$A$ leave a cone~$K$ invariant and let~$\rho = \rho(A)$. 
Then there exists a face~$\bar L$ of~$K$  invariant with respect to~$A$,  such that ${\rm dim}\, \bar L \le d-1$ and 
for every~$\bx \in K \setminus \bar L$, we have 
$\|A^j\bx\|\,  \dot{\asymp}\,  \rho^{\, j}, \ j\in \n$.  
\end{prop}
{\tt Proof.} If all eigenvalues of~$A$ have equal absolute values~$\rho$, 
then the assertion is obvious true with~$\bar L = \{0\}$. 
Otherwise, denote by~$\rho_2$ the second largest modulus of eigenvalues of~$A$. Thus, 
$\rho_2 < \rho$ and all 
vectors $\bx$ such that 
$\|A^j\bx\| \, \dot{\le}\, \rho_2^{\,j}$ form a proper invariant subspace~$V$ of
the operator~$A$. For every~$\bx\in \re^d$, we have $\|A^j\bx\| \, \dot{\asymp}\, \rho^{\, j}$, unless~$\bx \in V$. Denote $K' = K \setminus \{\bO\}$. 
If~$V$ does not intersect~$K'$, then the proposition follows 
with~$\bar L = \{\bO\}$. If it intersects the interior of~$K$, then we 
take an arbitrary point~$\bz \in {\rm int}\,K \, \cap \, V$
and consider the vector~$\bz - t\bv$, where $\bv$ is the Krein-Rutman eigenvector, 
for which $A\bv = \rho\, \bv$, and $t>0$ is small enough, so~$\bz - t\bv \in K$. 
We have $A^j(\bz - t\bv) \in K$ for all~$j$. 
On the other hand, $A^j(\bz - t\bv) \, = \, A^j\bz \, - \, tA^j\bv \, = \, 
A^j\bz \, - \, t\, \rho^{\,j}\bv\, = \, - t\, \rho^{\,j}\bv \, + \, o(\rho^{j})$
as~$j \to +\infty$.
However, $- t\, \rho^{\,j}\bv \, + \, o(\rho^{j}) \, = \, 
\rho^{\,j}\bigl( - t\bv \, + \, o(1)\bigr) \notin K$ for large~$j \to +\infty$, which is a contradiction. 
Thus, $V$ cannot intersect the interior of~$K$. It remains to consider the last case: 
$V$ intersects $K'$ but $V \cap {\rm int}\,K = \emptyset$. 
  By the convex separation theorem, $V$ can be separated from~$K$
  by a hyperplane~$H$. The set $H \cap K$ is a proper face 
  of~$K$ that contains the set~$V \cap K$. Since the operator~$A$
  respects the set~$V \cap K$, it respects its minimal 
  face, which will be denoted by~$\bar L$. Thus, $\bar L$
  is an invariant face, every  point $\bx\in K \setminus \bar L$
   does not lie in~$V$, 
  and hence~$\|A^j\bx\|\,  \dot{\asymp}\,  \rho^{\, j}$.

{\hfill $\Box$}
\medskip

\begin{cor}\label{c.20}
Suppose  a linear operator~$A$ leaves a cone~$K \subset \re^d$ invariant; 
then, for every~$\bx \in K$, we have 
\begin{equation}\label{eq.c20}
\bigl\|A^j\bx\bigr\|\quad \asymp\quad 
\bigl\|A^j\bigl|_{L}\bigr\| \quad \dot{\asymp}\quad  \rho^{\,j}\bigl(A\bigl|_{L}\bigr), \qquad j\in \z, 
\end{equation}
where  $L$  is the minimal invariant face of~$\bx$. 
\end{cor}
{\tt Proof.} It suffices to consider the case~$\|\bx\| = 1$ and~$L = K$, otherwise we 
restrict~$A$ to the corresponding face plane. If $\rho(A) = 0$, then $A$ is a nilpotent 
operator, $A^j = 0$ for $j\ge d$, and~(\ref{eq.c20}) holds.  
If~$\rho(A) > 0$, then by Proposition~\ref{p.20}, the assertion~(\ref{eq.c20})  holds 
for all~$\bx \in K$ which do not belong to a special invariant face of~$K$. 
However, by the assumption, $\bx$ does not belong to proper invariant faces,
which concludes the proof.

{\hfill $\Box$}
\medskip 

We are ready to prove the main result of this section, which gives us the key to  computing the regularity of refinable functions.

\begin{theorem}\label{th.45}
Let a linear operator~$A$ leave a cone~$K \subset \re^d$ invariant and 
 $\Gamma$ be an arbitrary subset of~$K$. Let also~$L$ be the minimal 
  invariant face for~$\Gamma$. Then for every~$\bx \in \Gamma$, 
 we have $\ \|A^j\bx\bigr\|\ \dot{\le}\  
 \rho^{\,j}(A|_{L}), \ j\in \n$. Moreover,  
there exists~$\bx_0 \in \Gamma$ and~$C>0$ such that $\ \|A^j\bx_0\bigr\|\, \ge \,   
 C\, \rho^{\, j}(A|_{L}), \ j\in \n$. 
\end{theorem}
{\tt Proof.} It suffices to consider the case~$L = K$, otherwise we 
restrict~$A$ to the corresponding face plane. Denote~$\rho = \rho(A)$
and take the invariant face~$\bar L$  from  
Proposition~\ref{p.20}.  Since $\Gamma$ does not lie on a proper invariant face
of~$K$, 
it follows that there exists~$\bx_0 \, \in \, \Gamma \setminus \bar L$, $C > 0$, for 
which~$\, \|A^j\bx_0\bigr\|\, \ge \, \quad  
 C\, \rho^{\, j}(A), \ j\in \n$. On the other hand, 
 $\, \|A^j\bx\bigr\|\, \dot{\le} \, 
 \rho^{\, j}$ for all~$\bx \in L$, which completes the proof.

{\hfill $\Box$}
\medskip

\bigskip 

\section{Cones of nonnegative  polynomials} \label{section:pos-poly}

\bigskip

Consider a finite subset~$\Omega \subset \z^n$ symmetric about the origin
and containing zero. 
Denote by~$\cP_{\Omega}$ the space of real-valued trigonometric 
polynomials on~$\re^n$ with the spectrum on~$\Omega$: 
$$
\cP_{\Omega} \quad =  \quad \left\{ \ \bp(\bxi)\ = \ \sum_{\bk \in \Omega}\, p_{\bk} \, e^{2\pi i \, (\bk, \bxi)}\quad :  
\quad  
  p_{-\bk}\, = \,  p_{\bk}\, \in \re\, , \ \bk \in \Omega\ \right\}.
$$ 
We usually assume that~$\Omega$ is fixed and use the simple notation 
~$\cP = \cP_{\Omega}$. 
We have~${\rm dim}\, \cP\ =\ \bigl(1+|\Omega|\bigr)/2$.  
Consider the set~$\cP_+$ of  nonnegative polynomials from~$\cP$, i.e., 
$\bp(\bxi) \ge 0, \ \bxi \in \re^n$. For a given closed set~$S \subset \re^n$, 
we denote 
$$
\cP(S) \quad = \quad \Bigl\{\, \bp \in \cP_{\Omega} \ : \ \bp(\bxi) = 0, \  \bxi \in S\ \Bigr\}\, , 
\qquad \cP_+(S) \, = \, \cP(S)\cap \cP_+\, . 
$$ 
The same notation is used for compact sets on the period: $\, S\subset \Tn$. 
Clearly, the sets $\cP(S)$ are linear subspaces of~$\cP$ and $\cP_+(S)$ are cones
in those subspaces. 
\begin{prop}\label{p.60}
The set~$\cP_+$ is a closed solid pointed cone in~$\cP$.
Every face~$L$ of~$\cP_+$ has the form~$L = \cP_+(S)$ for a suitable 
 compact subset~$S\subset \Tn$. 
\end{prop}
{\tt Proof.} The closeness and pointedness of~$\cP_+$ is obvious.  
To prove that this cone is solid, one can note that for every~$\bq \in \cP$,
such that~$\|\bq\|_{C(\Tn)} < 1$, we have $1 + \bq \in \cP_+$, hence~$\cP_+$ contains an open set. Now observe, that~$\cP_+$ is an intersection of 
the half-spaces~$\{\bp \in \cP \ | \ \bp(\bxi) \ge 0\}$ over all~$\bxi \in \Tn$. 
Therefore, every plane of support of the cone~$\cP_+$ has the 
form~$\cP(\{\bxi\})\, = \, \bigl\{\, \bp \in \cP\ : \  \bp(\bxi) = 0\bigr\}$ and hence,  every face plane of~$\cP_+$ is an intersection of some of those planes. 
Thus, every face of the cone~$\cP_+$ has the form~$\cP_+(S)$ for certain~$S \subset \Tn$.  

{\hfill $\Box$}
\medskip 

The proof of the following fact is the same as for~Proposition~\ref{p.60}: 
\begin{prop}\label{p.70}
For every compact set~$S\subset \Tn$, all faces 
of the cone~$\cP_{+}(S)$ are the sets $\cP_+(S'), \ S'\subset S$. 
\end{prop}

For a given~$\ell \ge 0$, we denote by~$\cP^{(\ell)}$ 
the set of polynomials $\bp \in \cP$ that have zero of order at least~$2(\ell+1)$ at the origin, i.e.,  
$\ \bp(\bO) = \bp'(\bO) = \cdots = \bp^{(2\ell + 1)}(\bO) = 0$. The equality $\bp^{(k)}(\bx) = 0$ means that all partial derivatives of order $k$ at the point~$\bx$ are equal to zero.  

We consider zeros of even order only, because nonnegative 
polynomials can have only such zeros at the origin.  The cone of nonnegative polynomials~$\cP^{(\ell)}_+ \, = \, \cP_+ \cap \cP^{(\ell)}$ is 
full-dimensional in~$\cP^{(\ell)}$. Similarly we define~$\cP^{(\ell)}_+(S)$
and prove the analogues of Propositions~\ref{p.60} and~\ref{p.70}  
for the space~$\cP^{(\ell)}$. We will refer to those propositions mentioning that 
we use them for~$\cP^{(\ell)}$ instead of~$\cP$. Sometimes we use 
the short notation~$P_{\ell} = P^{(\ell)}$.

\section{Regularity of 
 refinable functions with an arbitrary dilation matrix}\label{section:regularity}

\begin{center}
\textbf{5.1 Preliminary facts and notation}

\end{center}

We deal with the refinement equation~(\ref{eq.ref0}) 
with complex coefficients $\{c_{\bk}\}_{\bk \in \z^n}$ 
supported on a finite set~$Q\subset \z^n$ and such that
 $\sum_{\bk \in Q} c_{\bk} \, = \, m$. 
  For the sake of simplicity, we restrict ourselves to 
 the case or real coefficients, although all our results are 
 directly extended to the complex case. The unique compactly supported 
 solution $\varphi \in \cS_0'$ (refinable function) is normalized 
 by the condition~$\int_{\re^n}\varphi \, d\bx = 1$.   
  The equation   {\em satisfies the sum rules of order}~$\ell \ge 0$ if 
the $m-1$ points $\bxi = M^{-T}\bd_k^*, \ \bd_k^* \in D^*\setminus  \{\bO\}$ are  roots of order at least $\ell+1$  of the mask~$\bc(\bxi)$.  Recall that $D$ and $D^*$ are the sets of digits for the 
matrices~$M, M^T$ respectively.

\smallskip 

We make the standard assumption that 
the refinable function~$\varphi$ is {\em stable} i.e., its integer shifts~$\{\varphi(\vardot - \bk)\}_{\bk \in \z^n}$ are linearly independent. Equivalently, those 
shifts form a Riesz basis of their linear span. Without loss of generality we set~$\bd_0 = \bd_0^* = \bO$. 
It is known that if the solution $\varphi$ is stable and belongs to $W^{\ell}_2(\re^n)$, then the refinement equation 
satisfies the sum rules of order~$\ell$. In particular, if $\varphi \in L_2(\re^n)$, 
then the sum rules of order zero are fulfilled, i.e., $\bc\, \bigl( M^{-T}\bd_k^*\bigr) \, = \, 0, \ k=1, \ldots , m-1$.  
In terms of the coefficients, this condition reads: 
\begin{equation}\label{eq.sumrules0}
\sum_{\bj \in \z^n}\, c_{M \bj - \bd_k} \ = \ 1 \, , \qquad k = 0, \ldots , m-1\, . 
\end{equation}  
It will always be  assumed that the sum rules of order zero are satisfied. 
Now define several special subsets of~$\re^n$ and of~$\z^n$. 
Let 
$$
Y\quad =\quad \left\{ \  \sum_{j =1}^{\infty} M^{-j}\, \bs_j \ : \quad \bs_j \in Q\, , \ j \in \n \ \right\}\, .
$$ 
This is a compact self-similar set generated by the  affine operators~$A_{\bs}\bx \, = \, M^{-1}(\bx + \bs), \, \bs \in Q$. This means that~$Y$ 
is a unique fixed point of the set-valued map~$AX \, = \, \cup_{\bs \in Q} A_{\bs}X$, 
i.e., an affine fractal.  The transition operator~$\cF$
respects the space of functions supported on the set~$Y$, therefore, 
${\rm supp}\, \varphi \, \subset Y$, see~\cite{CHM1, CDM} for details. Let $\Omega \, = \, \bigl(Y - Y \bigr)\cap \z^n$,  
where $A-B \, = \, \{\ba-\bb \, : \ \ba\in A, \bb \in B\}$ is the standard 
Minkowski subtractions of two sets. Thus, $\Omega$ is a symmetric finite (since~$Y$ is compact)
subset of~$\z^n$ containing zero.  This set is uniquely defined by the support of the mask~$Q$
and can be constructed by an algorithm presented in Appendix.  
It follows from the definition that $\Omega \, = \, \{\bk \in \z^n \, : \ Y \cap ( Y + \bk)\, \ne \, \emptyset\}$. Therefore,~$\Omega$ contains all integer points from the 
support of the function~$\varphi * (\varphi\_)$ (the symmetrized autoconvolution 
of~$\varphi$), where $\varphi\_(\vardot) = \varphi (-\vardot)$.

We consider the linear operator~$\cT$ defined in~$L_2(\Tn)$ by formula~(\ref{eq.T0}). 
This operator maps the space of polynomials~$\cP_{\Omega}$ to itself~\cite{D, E, J} and 
in the sequel we restrict~$\cT$ to this subspace. 
If the mask~$\bc$  satisfies the sum rules of order~$\ell$, then $\cT$  respects also all the subspaces~$\cP^{(k)}, \, k \, = \, 0, 1, \ldots , \ell$, see~\cite{H3, JZ}.  Clearly, it  
 respects the positive cones~$\cP^{(k)}_+\, $ as well.  
For an arbitrary compactly supported distribution~$f\in \cS'_0$, 
the  periodization of the square 
of its Fourier transform $\Phi_f$ defined in~(\ref{eq.Phi0}) is a non-negative trigonometric polynomial
with coefficients~$p_{\bk} \, = \, 
[\, f* f\_\, ](-\bk)$. 
So, if ${\rm supp}\, f\, \subset Y$, then $\Phi_f\, \in \, \cP_{+}$.
For every~$f \in \cS_0'$, we have~$\cT (\Phi_f) \, = \, \Phi_{\cF f}$.

\medskip

For each~$\bp \in \cP_+$, we are interested in 
the asymptotic behaviour of the sequence~$\{\cT^k\bp\}_{k\in \n}$
 as $k \to \infty$. 
In view of Theorem~\ref{th.45}, it is defined by the spectral radius of the 
operator~$\cT$ restricted to the minimal invariant 
face of the cone~$\cP_+$ that contains~$\bp$. That face, in turn, is determined by
the set of  zeros~$S$ of the polynomial~$\bp$  
(Proposition~\ref{p.60}). This set~$S$ must possess a special structure because 
the space of polynomials that vanish on~$S$ has to be  invariant with respect to the 
operator~$\cT$.   
We need to reveal this structure and find this set.

For arbitrary~$s=1, \ldots , q$, 
we consider the set 
$\,  \pi \, \bigl\{ \nill( \widehat\varphi) \, , \,  J_s^{\perp}\bigr\}\,$
that consists of points~$\bxi \in \re^n$
such that $\bxi + \bk  \ \in \ \nill( \widehat\varphi) \, \cup \,   J_s^{\perp}$
for all~$\bk \in \z^d$. Let us remember that $\nill(f)$ denotes the set of zeros of an analytic function~$f$ and $J^{\perp}$ is an orthogonal complement to 
the subspace~$J$ in~$\re^n$. 
 Thus, the periodized point~$\bxi$ is contained in the union of 
the sets $\nill( \widehat\varphi)$ and $J_s^{\perp}$.  
The stability of the function~$\varphi$ implies that its Fourier transform~$\widehat\varphi$
does not have periodic zeros. This yields that~$\pi \, \bigl\{ \nill( \widehat\varphi)\bigr\}\, = \, \emptyset$. So, only the hyperplane~$J_s^{\perp}$ makes this set 
nonempty. 
Indeed, the set~$\nill( \widehat\varphi)$ contains all integer points except for zero while 
$J_s^{\perp}$ contains zero, consequently,~$\bO \in \pi \, \bigl\{ \nill( \widehat\varphi) \, , \,  J_s^{\perp}\bigr\}$.

Denote by~$P_{\ell ,s}$   
the space of polynomials from~$\cP_{\Omega}$    
that have zeros of order at least~$2(\ell + 1)$ at all  points of  
the set  $\pi \, \bigl\{ \nill_{\ell}( \widehat\varphi) \, , \,  J_s^{\perp}\bigr\}\,$ (see subsection~1.4 for the notation). The cone of nonnegative polynomials in~$P_{\ell ,s}$ is referred to as~$K_{\ell , s}$. 
 This is a face of the cone~$\cP_+$ with the corresponding face plane~$P_{\ell ,s}$.

\begin{lemma}\label{l.50}
If the refinement equation satisfies sum rules of order~$\ell$, 
then for every~$s = 1, \ldots , q$, the  space~$P_{\ell, s}$ is an invariant 
subspace of the operator~$\cT$ and $K_{\ell, s}$ is its invariant cone in this subspace. 
\end{lemma}  
{\tt Proof} is in Appendix.  
\medskip 

\begin{prop}\label{p.90}
Let a refinement equation satisfy the sum rules of order~$\ell$. 
Then, for each~$s = 1, \ldots , q$, there is an alternative: 
either~$P_{\ell, s} = \cP^{(\ell)}$ or 
there exists a non-integer vector~$\bz \in \re^n$ orthogonal to~$J_s$ such that 
for every~$\bk \in \z^n$, the vector~$\bz+\bk$
is orthogonal to~$J_s$ or is a zero 
of the function  $\widehat \varphi$ 
of order at least~$\ell+1$. 
\end{prop}  
{\tt Proof.} As we already noted, the set
$\pi \, \bigl\{ \nill_{\ell}( \widehat\varphi) \, , \,  J_s^{\perp}\bigr\}\, $  always  contains~$\{\bO\}$. 
On the other hand, all polynomials from $\cP^{(\ell)}$
have a zero of order~$\, \ge\, 2(\ell + 1)$ at the point~$\bO$. 
Hence,~$P_{\ell, s} \, \ne  \, \cP^{(\ell)}$ if and only if the set 
$\pi \, \bigl\{ \nill_{\ell}( \widehat\varphi) \, , \,  J_s^{\perp}\bigr\}\, $
contains a non-integer point~$\bz$. Since the function~$\varphi$ is stable, it follows that $\pi \, \bigl\{ \nill_{\ell}( \widehat\varphi)\bigr\}\, = \, 
\emptyset$. Therefore, $\bz \in J_s^{\perp}$, which concludes the proof.

{\hfill $\Box$}
\medskip 

\bigskip

\begin{center}
\textbf{5.2. The $L_2$-regularity}

\end{center}

\bigskip

Now we prove the formula for the~$L_2$-regularity and in the next section attack the 
higher order regularity. The main idea is to work separately with 
regularities along 
spectral subspaces~$J_1, \ldots , J_q$ of~$M$. Then the regularity of~$\varphi$ is equal to the minimal regularity over these subspaces.

We denote by~$\alpha_0(f)$ the H\"older exponent 
(or the {\em H\"older exponent of zero order}) of a function~$f \in L_2(\re^n)$. 
The index zero is added to distinguish from the    H\"older exponent of higher order.  

For a given nontrivial subspace~$J\subset \re^n$ (the case~$J=\re^n$ is allowed), 
we define the H\"older exponent in~$L_2$ of the  function~$f \in L_2(\re^n)$
along~$J$: 
$$
\alpha_{\,0,J} \quad = \quad \sup\  
\Bigl\{\alpha \ge 0 \ : \ \|f(\vardot + \bh) - f(\vardot)\| \ \le \ C\, \|\bh\|^{\, \alpha}
\, ,  \  \bh \in J \Bigr\}\, .
$$
 For the subspaces~$J_s$, we use the short 
notation~$\alpha_{\,0,s}\, = \, \alpha_{\,0,J_s}\, , \ s=1, \ldots , q$.

 The spectral radius of the operator~$\cT$ restricted to its invariant 
 subspace~$P_{0, s}$ is denoted by~$\rho_{\,0,s}$. 

\begin{theorem}\label{th.50}
For an arbitrary refinable function~$\varphi$ with the dilation matrix~$M$, we have 
\begin{equation}\label{eq.Hold0s}
\alpha_{\,0,s}(\varphi)\quad = \quad \frac12\, \log_{1/r_s} \, \rho_{\,0,s}\, , \qquad 
s=1, \ldots , q, 
\end{equation}
where $J_s$ are the spectral subspaces of~$M$. For the H\"older exponent
in~$L_2(\re^n)$, we have 
\begin{equation}\label{eq.Hold0}
\alpha_0 (\varphi)\quad  = \quad  \frac12\,  \min_{s=1, \ldots , q}\log_{\, 1/r_s} \, \rho_{\,0,s}
\end{equation}
\end{theorem} 
Note that~$P_{0,1}$ can be a proper subspace of~$\cP^{(0)}$
only if the mask possesses a special set of zeros. 
That is why, 
for most of refinement equations, we have~$P_{0,1} = 
\cP^{(0)}$, and hence, the spectral radius~$\rho_{\,0,1}$ 
coincides with~$\rho_0 \, = \, \rho \, \bigl(\cT|_{\cP^{(0)}} \bigr)$. 
Since $r=r_1$ is the maximal spectral radius among all~$r_s, \, s =1, \ldots , q$, 
we see that the minimum in the formula~(\ref{eq.Hold0}) is attained for~$s=1$. 
Therefore,  
$\alpha_0 (\varphi) \, = \, \frac12\, \log_{\, 1/r_1} \, \rho_{0}$. 
Applying Proposition~\ref{p.90} for the subspace~$J_1$, we obtain 
\begin{theorem}\label{th.55}
The equality $\alpha_0 (\varphi)\, = \, \frac12\,  \log_{1/r} \, \rho_0$ holds,  
unless there exists a non-integer point~$\bz \in J_1^{\perp}$ 
such that, for every integer shift~$\bz + \bk, \, \bk \in \z^n$, we have 
either~$\bz + \bk \, \in \, J_1^{\perp}$ or  $\widehat \varphi (\bz + \bk) = 0$. 
\end{theorem}
There are, however,  many examples when $P_{0,1}$ is a proper subspace 
of~$\cP^{(0)}$ and, respectively, $\rho_{\,0,1} < \rho_0$. In this case,  
the H\"older exponent can be larger than~$\frac12\, \log_{\, 1/r_1} \, \rho_{0}$
and is equal to the minimal value~$\, \frac12\, \log_{\, 1/r_s} \, \rho_{\,0, s}$
among~$s\, =\, 2, \ldots , q$.

The proof of Theorem~\ref{th.55}   requires several further notation and auxiliary results. 
In the sequel of this section we use the simplified notation~$P = P_{0}$
(the subspace of polynomials from~$\cP_{\Omega}$ such that~$\bp(0) = \bp'(0) = 0$),
$P_s = P_{0,s}$ (the space of polynomials from~$P$ that vanish on the set
$\pi \{\nill (\widehat \varphi), J_s^{\perp}\}$). 
We also denote $K = K_0$ and $K_s = K_{0,s}$, those  are  the cones of nonnegative polynomials 
in~$P$ and in~$P_s$ respectively. 
\smallskip

An $\varepsilon$-extension~$Y_{\varepsilon}$ of~$Y$ is a compact set 
containing a neighbourhood of~$Y$ such that~$Y_{\varepsilon} \subset Y + B(0, \varepsilon)$
and~$AY_{\varepsilon} \subset Y_{\varepsilon}$, where~$AY_{\varepsilon} \, = \, \cup_{\bs \in Q} A_{\bs} Y_{\varepsilon}$.  Let us remember that~$AY = Y$. For every 
$\varepsilon$-extension~$Y_{\varepsilon}$,  the transition operator~$\cF$
respects the space of functions supported on the set~$Y_{\varepsilon}$. 
\begin{lemma}\label{l.60}
For every~$\varepsilon> 0$,  there exists an $\varepsilon$-extension of~$Y$. 
\end{lemma}
{\tt Proof.} Take a number~$t>0$ such that the diameter $tY $ is less than~$\varepsilon$
and define~$Y_{\varepsilon} = Y+tY$ (the Minkowski sum). This set contains 
a neighbourhood of~$Y$ because~$tY$ has a nonempty interior and, clearly, 
$AY_{\varepsilon} \subset Y_{\varepsilon}$. 

{\hfill $\Box$}
\medskip 

\begin{remark}\label{r.80}
{\em If $Y$ is not convex, then $Y+tY$ is not necessarily  equal to~$(1+t)Y$. 

}
\end{remark}

Let us remember that for every~$f \in \cS_0'$, the function~$\Phi_f$
is a trigonometric polynomial and if~${\rm supp}\, f\, \subset \, Y$, then 
$\Phi_f \in \cP_{\Omega}$, where $\Omega = (Y-Y)\cap \z^n$.  Since $Y$ is compact, it follows that 
a sufficiently small enlargement of the set~$Y$ does not change the set~$\Omega$. 
Thus, the set~$(Y_{\varepsilon}-Y_{\varepsilon})\cap \z^n$ 
coincides with~$\Omega$, whenever $\varepsilon$ is small enough.  
Fix such~$\varepsilon$ and an $\varepsilon$-extension~$Y_{\varepsilon}$. 
For each~$s=1, \ldots , q$, consider the set of differences 
of the function~$\varphi$ under small shifts along~$J_s$: 
\begin{equation}\label{Us}
U_s \ =  \ \Bigl\{f(\cdot )\, = \,  \varphi(\vardot + \bh)\, - \, \varphi(\vardot )\ :
\quad \bh \in J_s \, , \ {\rm supp}\, f\, \subset \, Y_{\varepsilon}\, \Bigr\}\, . 
\end{equation}
Those sets  are invariant with respect to the transition operator~$\cF$. 
Observe  that $U_s$ are not linear subspaces. 
Their linear spans 
may have nontrivial intersections, may  even all  coincide, while, recall,  
the spaces~$J_s, \, s=1, \ldots , q$, form a decomposition of~$\re^n$ into a direct sum. 
Every~$f\in U_s$
is defined by the translation vector~$\bh \in J_s$, that is why we sometimes 
denote~$f=f_{\bh}$. 

Now we are ready to prove  Theorem~\ref{th.50}. 
\medskip

{\tt Proof of Theorem~\ref{th.50}}. Consider the space~$J_s$ for some~$s = 1, \ldots , q$. 
For an arbitrary function~$f\in U_s$, the Plancherel theorem and the equality
$\Phi_{\, \cF f}\, = \, \cT \, \Phi_f, \, f\in \cS_0'$, imply
\begin{equation}\label{eq.trans3}
\bigl\|\cF^{j}\, f\bigr\|_2^2 \quad = \quad  \bigl\|\widehat {\cF^{j}\, f}\bigr\|_2^2 \quad = \quad 
\bigl\|\cT^j \Phi_f\, \bigr\|_1\, .   
\end{equation} 
Moreover,  the function~$\Phi_f$ is a polynomial from~$\cP_+$. 
Let $V_s$ be the minimal face of the cone~$\cP_+$ containing 
the set 
$\Gamma = \{ \Phi_f\, :\,  \ f  \in U_s\}$. Since~$\cF U_s \subset U_s$, it follows that 
$\cT V_s \subset V_s$.   
Let us show that~$V_s = K_s$. For each~$f\in U_s$, the Fourier transform~$\widehat f(\bxi)\, = \, 
\bigl( e^{2\pi i (\bh, \bxi)}\, - \, 1\bigr)\widehat{\varphi}(\bxi)$ vanishes on the family of parallel 
hyperplanes~$\{\bxi \in \re^n \, | \, 
(\bh, \bxi) \in \z \}$.  
This property holds for all sufficiently small~$\bh \in J_s$ precisely 
when~$f$ vanishes on the set of all sufficiently small vectors~$\bxi \in \re^n$ 
orthogonal to the subspace~$J_s$. This means that it vanishes 
on the orthogonal 
complement of~$J_s$ in~$\re^n$. Each~$\Phi_f \in \Gamma$ generating 
$V_s$ 
  vanishes on the union of the sets~$\nill (\widehat{\varphi})$
and~$\bh^{\perp}$. The intersection of all those sets over all 
feasible~$\, \bh\in J_s$
is $\nill\, (\widehat{\varphi})\, \cup \, J^{\perp}_s$. By Proposition~\ref{p.60}, 
every face of~$\cP_+$ consists of nonnegative polynomials from~$\cP_{\Omega}$
that vanish on a certain closed set~$S \subset \re^n$. Therefore, the face~$V_s$ is 
the set of polynomials from~$\cP_+$ that vanish on the set~$\pi \, \{\nill(\widehat{\varphi}), J_s^{\perp}\}$, i.e., $V_s = K_s$.

Since all eigenvalues of the matrix~$M$
restricted to the subspace $J_s$ are equal to~$r_s$, 
it follows that there exists a constant~$c_1> 0$ such that
$$
c_{1\, }r_s^j\ \|\bh\| \quad \le \quad \|M^j\bh\| \quad \dot{\le} \quad r_s^j\, \|\bh\|\, , \qquad \bh \in J_s, \ j \, \in \, \n\, . 
$$
Hence, for every~$\delta > 0$, there exists a constant~$c_2$
such that~$\|M^j\bh\| \ \le \ c_2\, (r_s+\delta)^j\, \|\bh\|, \, j\in \n$.  
Further, choose some~$\tau < \varepsilon$
such that the $\tau$-neighbourhood of~$Y$ is contained in~$Y_{\varepsilon}$. 
Take arbitrary~$\bh \in J_s, \, \|\bh\| < \tau$, and the smallest
 integer $j$ such that 
$\|M^{j+1}\bh\| \, \ge  \, \tau$. Since $\|M^{j}\bh\| \, <  \, \tau$, 
it follows that the function~$f=\varphi(\vardot + M^j\bh) -  \varphi(\vardot)$ belongs 
to~$U_s$. Applying~(\ref{eq.trans3}), we obtain 
$$
\Bigl\|\varphi(\vardot + \bh) \, - \, \varphi(\vardot)\Bigr\|^2 \ = \ 
\Bigl\|\cF^{j}\, \bigl( \varphi(\vardot + M^j\bh) -  \varphi(\vardot)\bigr) \Bigr\|^2
\ = \ \Bigl\|\cF^{j}\, f \Bigr\|^2 \ = \ \Bigl\|\cT^j \Phi_f\, \Bigr\|_1\, . 
$$
On the other hand,~$K_s$ is the minimal invariant face of~$\cP_+$ that contains $\Gamma$. 
 By Theorem~\ref{th.45}, for every $\Phi_f\in \Gamma$, 
we have~$\bigl\|\cT^j \Phi_f\, \bigr\|_1 \, \dot{\le}\, (\rho_{\, 0,s})^{\, j}\|\Phi_f\|_1$.  Furthermore, 
$\|\Phi_f\|_1 \, = \, \|\widehat f\|^2 \, = \, \| f\|^2 \, \le \, 
4\|\varphi\|^2$, because $\|f\|\, = \, \| \varphi(\vardot + M^j\bh) -  \varphi(\vardot)
\|\, \le \, 2\, \|\varphi\|.$ Thus, 
\begin{equation}\label{eq.trans9}
\bigl\|\varphi(\vardot + \bh) \, - \, \varphi(\vardot)\bigr\|^2 \quad \dot{\le} \quad 
 (\rho_{\, 0,s})^{\, j}\, . 
\end{equation} 
On the other hand, $\ \tau \, \le \, \|M^{j+1}\bh\| \ \le \ 
c_2(r_s + \delta)^{\, j+1}\|\bh\|$. 
Therefore,  $\|\bh\| \, \ge \, \frac{\tau}{c_{2}}\, (r_s+\delta)^{j+1}$. Combining this 
with~(\ref{eq.trans9}) and taking the limit as~$j\to \infty$, we obtain 
$\, \alpha_{\, 0, s} \, \ge \, \frac12\, \log_{1/(r_s+\delta)}\, \rho_{\, 0,s}$. This is true for every~$\delta > 0$, 
consequently,~$\, \alpha_{\, 0, s} \, \ge \, \frac12\, \log_{1/r_s}\, \rho_{\, 0,s}$. 

To establish the inverse inequality we 
invoke Theorem~\ref{th.45}: 
there exists~$\Phi_f\in \Gamma$ 
such that~$\|\cT^j \Phi_{f}\| \ \ge  \ C\, (\rho_{\,0,s})^j$. 
Applying~(\ref{eq.trans3}), we get 
 $\, \bigl\|\cF^{j}\, f\bigr\|^2 \, = \,  \|\cT^j \Phi_{f}\|    \, \ge \, C\, (\rho_{\,0,s})^j$. 
Since $f(\vardot) \, = \, \varphi(\vardot + \bh) - \varphi(\vardot)$ for some 
$\bh \in J_s$, we see that 
$$
\bigl\| \varphi(\vardot + M^{-j}\bh) - \varphi(\vardot) \bigr\| \ = \ 
\bigl\| \cF^j ( \varphi(\vardot +\bh) - \varphi(\vardot) ) \bigr\|\ = \ 
\bigl\| \cF^j f \bigr\|\ \ge \ C\, (\rho_{\,0,s})^{j/2}\, . 
$$
 On the other hand, $\, \|M^{-j}\bh\| \, \le\, c_{1}\, (r_s-\delta)^{-j}\|\bh\|$, which implies that~$\alpha_{\, 0, s} \, \le \, 
\frac12\, \log_{1/(r_s-2\delta)}\, \, \rho_{\, 0,s}$.
 Taking the limit as~$\delta \to 0$, we  conclude the proof.

{\hfill $\Box$}
\medskip

\begin{center}
\textbf{5.3 The higher order regularity}

\end{center}

For an integer $k \ge 0$ and for a given nontrivial subspace~$J\subset \re^n$, we 
denote by $W^{k}_{2}(\re^n, J)$ the space of functions
$k$ times differentiable along~$J$. This means that  for every 
multivector  $\bar \bb \, = \, (\bb_1, \ldots , \bb_{k}) \in  J^{\,k}$, 
where $J^{\,k} \, = \, J\times \cdots \times J$ ($k$ multipliers),  
$\bb_i \ne 0, \, i = 1, \ldots , k$, the mixed partial derivative $f_{\bar \bb}\, = \, \frac{\partial}{\partial \bb_1 \cdots \partial \bb_{k}}\, f$ along the vectors $\bb_1,  \ldots , \bb_{k}$  belongs to~$L_2(\re^n)$.  
For every  $f\in W^{k}_{2}(\re^n, J)$, the {\em H\"older exponent of order~$k$ along~$J$}
is defined as follows:  
$$
\alpha_{k, J} \quad = \quad k\ +\  \sup\  
\Bigl\{\alpha \ge 0 \ : \ \|f_{\bar \bb}(\vardot + \bh) - f_{\bar \bb}(\vardot)\| \ \le \ C\, \|\bh\|^{\, \alpha}
\, ,  \  \bar \bb \in J^{\, k} \Bigr\}\, . 
$$
For the subspaces~$\{J_s\}_{s=1}^q$ (the spectral subspaces 
of the dilation matrix~$M$), we use the short 
notation~$\ \alpha_{k, J_s}\, = \, \alpha_{\, k, s}\, , \ s=1, \ldots , q$.

If $J=\re^n$, then we obtain the standard definition of the 
{\em H\"older  exponent of the $k$th order}~$\alpha_k(f)$
as the number~$k$ plus the minimal H\"older  exponent (of  order zero)
of all mixed derivatives of order~$k$.  
The (generalized) H\"older  exponent of a function~$f\in L_2(\re^n)$
  is $\alpha(f)\, =\, \alpha_{\,k}(f)$, 
where $k$ is the largest integer, for which $f\in W_2^k(\re^n)$.

From the representation by the direct sum~$\re^n \, = \, \oplus \sum_{s=1}^q J_s$, 
it follows that if  a function~$f$ belongs to~$W_2^k(\re^n, J_s)$ for each~$s$, then it belongs 
to~$W_2^k(\re^n)$ and the H\"older exponent $\alpha_k(f)$ is equal to 
$\min_{s=1, \ldots , q} \alpha_{\,k, s}$.  The  computation of each exponent of regularity~$\alpha_{\,k, s}$ will be done similarly to the 
previous section, in terms of spectral radii of the operator~$\cT$
on faces of the cone of nonnegative polynomials~$\cP_+$.

In the statement of the main result we invoke the notation 
of the space~$P_{k,s}$,  the cone  
$K_{k,s}$, and the spectral radius~as~$\rho_{k,s}$ of the operator~$\cT|_{P_{k,s}}$, 
see subsection~5.1. Recall that~$P_{k,s}$ is  respected by~$\cT$, provided 
the refinable function~$\varphi$ satisfies the sum rules of order~$k$
(Lemma~\ref{l.50}). We also use the short notation~$P_k = P^{(k)}$. 

\begin{theorem}\label{th.50h}
A refinable function~$\varphi$ satisfying the sum rules of order~$k$
belongs to~$W^{k}_2(\re^n, J_s)$
precisely when  $\frac12\, \log_{1/r_s} \rho_{k, s} \, >  \, k$. 
In this case $\alpha_{\,k,s}(\varphi) \, = \,  \frac12\, \log_{1/r_s} \rho_{k, s}$. 

The function~$\varphi$ belongs to~$W_2^{k}(\re^n)$ if and only 
it satisfies the sum rules of order~$k$ and 
 belongs to~$W^{k}_2(\re^n, J_s)$ for all~$s=1, \ldots , q$. 
The H\"older exponent of~$\varphi$ is equal to 
\begin{equation}\label{eq.Hold}
\alpha(\varphi) \ = \ \min_{s=1, \ldots , q} \, \frac12\, \log_{1/r_s} \rho_{k, s}\, , 
\end{equation}
where $k$ is the maximal number such that~$\varphi \in W^k_2(\re^n)$. 
\end{theorem} 
Applying Proposition~\ref{p.90} for the subspace~$J_1$, we obtain 
\begin{theorem}\label{th.55h}
The equality $\alpha (\varphi)\, = \, \frac12 \, \log_{1/r} \, \rho_{k}$,  
where~$r=\rho(M), \, \rho_k = \rho(\cT|_{P_k})$ and 
$k$ is the maximal number such that 
$\frac12 \, \log_{1/r} \, \rho_{k} \, > \, k$, holds, 
unless there exists a non-integer point~$\bz \in J_1^{\perp}$ 
such that every integer shift~$\bz + \bk \, , \ \bk \in \z^n$
is either on~$J_1^{\perp}$ or is a  zero of the function~$\widehat \varphi$
of order at least~$k+1$. 
\end{theorem}

The proof of Theorem~\ref{th.55h} is realized in the same way as the proof of Theorem~\ref{th.50}, but replacing 
the difference~$\varphi(\vardot + \bh) - \varphi(\cdot)$ 
by the higher order differences~$\Delta_{\, \bar \bh}\varphi$. We begin with necessary definitions.

For a given multivector~$\bar \bh = (\bh_0, \ldots , \bh_{k})\in (\re^n\setminus \{0\})^{k+1}$, 
we denote by~$\Delta_{\, \bar \bh}f (\vardot)$ the~$(k+1)$st order  difference 
along the vectors~$\bh_i$. It is  defined inductively as 
$$
\Delta_{\, (\bh_0, \ldots , \bh_i)}f(\vardot) \ = \ 
\Delta_{\, (\bh_0, \ldots , \bh_{i-1})}f(\vardot + \bh_i) \ - \ 
\Delta_{\, (\bh_0, \ldots , \bh_{i-1})}f(\vardot )\ , \qquad i=1, \ldots , k\, . 
$$
n what follows we always assume that all the vectors of shifts~$\bh_i$
are nonzero and denote~$\|\bar \bh\| \, = \, \max_{i=0, \ldots , k}\|\bh_i\|$. 
We need the following   well-known fact from the approximation theory~\cite{K, Mil}, see also~\cite[Theorem 4.6.14]{TB}. It is usually formulated for~$J = \re^n$, 
 but the proofs are the same  for proper subspaces~$J \subset \re^n$.   
\smallskip

\noindent \textbf{Theorem A}. {\em For an arbitrary function~$f\in W^k_2(\re^n, J)$, we have 
$$
\|\Delta_{\, \bar \bh}f\| \ = \ o(\|\bar \bh\|^{k}), \qquad  \bar \bh \in J^{k+1}, \  
\|\bar \bh\| \, \to\, 0\, .
$$ 
Conversely, if  $\|\Delta_{\, \tilde \bh}f\| \, = \, O\, \bigl( \|\tilde \bh\|^k\bigr), 
 \, \tilde  \bh \in J^{k}, \, \|\tilde  \bh\| \, \to\, 0$, then~$f$ 
 belongs to~$W^k_2(\re^n, J)$ after possible correction on a set of zero measure.  
 In this case 
\begin{equation}\label{alphak}
\alpha_{\,k,J}(f)\quad = \quad \sup\ \Bigl\{\, \alpha >0 \ : \ \|\Delta_{\, \bar \bh}f\| \, \le \, C\, \|\bar \bh\|^{\, \alpha}, \qquad \bar \bh \in J^{k+1} \, \Bigr\}\, . 
\end{equation}}

We use the $\varepsilon$-extensions~$Y_s$ of the set~$Y$ 
(Lemma~\ref{l.60}) and the following high-order analogues of the sets~$U_s$
from~(\ref{Us}): 
\begin{equation}\label{Ush}
U_{k, s} \ =  \ \Bigl\{ \, f = \Delta_{\, \bar \bh}\varphi\,  :
\quad \bar \bh \in J_s^{\, k+1} \, , \ {\rm supp}\, f\, \subset \, Y_{\varepsilon}\, \Bigr\}\, . 
\end{equation}
Thus, $U_{k, s}$  is the set of~$(k+1)$st order differences 
of the refinable function~$\varphi$ under small shifts along~$J_s$. 
All those sets are invariant with respect to the transition operator~$\cF$.

\smallskip

{\tt Proof of Theorem~\ref{th.50h}}.  
For arbitrary~$f\in U_{k,s}$,  the function~$\Phi_f$ is a polynomial from~$\cP_+$. 
We denote by $V_{k,s}$ the minimal face of the cone~$\cP_+$ containing 
the set 
$\Gamma \, = \, \{ \Phi_f\, :\,  \ f  \in U_{k, s}\}$. Since~$\cT \, \Gamma \subset \Gamma$, it follows that  
$\cT V_{k,s} \subset V_{k,s}$. To show that~$V_{k,s} = K_{k,s}$
we observe that the Fourier transform of every function~$f\in U_{k,s}$
has the form~$\widehat f(\bxi)\, = \, 
\widehat{\varphi}(\bxi)\prod_{a=0}^{k}\, \bigl( e^{2\pi i (\bh_a, \bxi)}\, - \, 1
\bigr)$ and hence,  vanishes on each family of parallel 
hyperplanes~$\{\bxi \in \re^n \, | \, 
(\bh_a, \bxi) \in \z \}, \ a = 0, \ldots , k$.  
This property holds for all sufficiently small~$\bh_a \in J_s$ precisely 
when~$f$ has zero of order at least~$k+1$ at each point of the subspace~$J_s^{\perp}$.  
Therefore, each polynomial~$\Phi_f$  has zeros of order at least $2(k+1)$ on 
 the union of the sets~$\nill_k (\widehat{\varphi})$
and~$\bh_a^{\perp}$. The intersection of all those sets over all 
feasible~$\bh_a\in J_s, \, a=0, \ldots , k$, 
is $\nill_k (\widehat{\varphi})\cup J^{\perp}_s$. By Proposition~\ref{p.60}, 
the face~$V_{k,s}$ coincides with  
the set of polynomials from~$\cP_+$ that have zeros of order 
at least~$2(k+1)$ on the set~$\pi \, \{\nill_k(\widehat{\varphi}), J_s^{\perp}\}$, which is $ K_{k,s}$. 

The remainder of the proof is the same as for Theorem~\ref{th.50}. 
Choose some~$\tau < \varepsilon$
such that the $(k+1)\tau$-neighbourhood of~$Y$ is contained in~$Y_{\varepsilon}$.  
Denote $M^j\bar \bh \, = \, 
\bigl(M^j\bh_0, \ldots ,  M^j\bh_k\bigr)$. For arbitrary~$\bar \bh \in J_s^{k+1}, \, \|\bar \bh\| < \tau$, we find  the biggest integer $j$ such that 
$\|M^{j+1}\bar \bh\| \, \ge \, \, \tau$. 
The function~$f\, =\, \Delta_{\, M^j \bar \bh}\, \varphi\, $ belongs to~$U_{k,s}$, and hence,  applying~(\ref{eq.trans3}), we obtain 
$$
\Bigl\|\Delta_{\, \bar \bh}\varphi\Bigr\|^2 \ = \ 
\Bigl\|\cF^{j}\, \Delta_{\, M^j \bar \bh}\varphi \Bigr\|^2
\ = \ \Bigl\|\cF^{j}\, f \Bigr\|^2 \ = \ \Bigl\|\cT^j \Phi_f\, \Bigr\|_1\, . 
$$
On the other hand, $\Phi_f$ belongs to the set~$\Gamma$, 
for which~$K_{k,s}$ is the minimal  invariant face. 
By Theorem~\ref{th.45}, we have~$\bigl\|\cT^j \Phi_f\, \bigr\|_1 \, \dot{\le}\, (\rho_{\, 0,s})^{\, j}\|\Phi_f\|_1$, and consequently $\bigl\|\Delta_{\, \bar \bh}\varphi \bigr\|^2 \, \dot{\le} \, 
 (\rho_{\, k,s})^{\, j}\, .$ 
Combining this with the inequality~$\|M^{j+1}\bar \bh\| \, \ge \, \, \tau$, and taking~$j\to \infty$, we conclude that 
$\, \alpha_{\, k, s} \, \ge \, \frac12\, \log_{1/r_s}\, \rho_{\, k,s}$. 

Further, by Theorem~\ref{th.45}, there exists~$\Phi_f\in \Gamma$ 
such that~$\|\cT^j \Phi_{f}\| \ \ge  \ C\, (\rho_{k,s})^j$. 
Applying~(\ref{eq.trans3}), we get 
 $\, \bigl\|\cF^{j}\, f\bigr\|^2 \, = \,  \|\cT^j \Phi_{f}\|    \, \ge \, C\, (\rho_{k,s})^j$. 
Since $f \, = \, \Delta_{\, \bar \bh}\varphi$ for some 
$\bar \bh \in J_s^{\, k+1}$, we see that 
$$
\bigl\| \Delta_{\, M^{-j} \bar \bh}\varphi \bigr\| \ = \ 
\bigl\| \cF^j  \Delta_{\, \bar \bh}\varphi \bigr\|\ = \ 
\bigl\| \cF^j f \bigr\|\ \ge \ C\, (\rho_{\,0,s})^{j/2}\, . 
$$
Combining this with the equality $\, \|M^j\bar \bh\| \, < \, \, \tau$ and taking $j\to \infty$, we conclude that~$\alpha_{\, k, s} \, \le \, 
\frac12\, \log_{1/(r_s+\delta)}\, \, \rho_{\, k,s}$.
Taking the limit as~$\delta \to 0$, we  complete the proof.

{\hfill $\Box$}
\medskip

\section{The tile B-splines}\label{section:splines}

As an example of application of  Theorem~\ref{th.50h}, we compute the regularity 
of the tile B-splines. We are going to show that, unless the dilation matrix
$M$ can be factored to a block-diagonal form by means of an integer basis, 
the H\"older exponent  is computed by a simple formula $\alpha (\varphi)\, = \, \frac12 \, \log_{1/r} \, \rho_{\ell}$, where~$\ell$ is the order of the B-spline. 
 Thus,  the exponent of regularity 
is computed without involving zeros from the set~$\pi\, \bigl\{\nill_\ell(\widehat \varphi)\, , \, 
J_s^{\perp}\, \bigr\}$. The proof is based on two arguments: geometrical (tile cylinders) and 
combinatorial (rational subspaces in~$\re^n$). First we introduce tile cylinders 
and prove their basic  properties. We show that if the assumptions 
of~Theorem~\ref{th.55h} are not satisfied, the tiling by translates of~$G$ contains a 
nontrivial cylinder parallel to an eigenspace of~$M$. Then we prove that 
this eigenspace is a rational subspace of~$\re^n$. This means that 
if $M$ does not have rational eigenspaces, then the assumptions of 
Theorem~\ref{th.55h} are satisfied.

\bigskip 

\begin{center}
\textbf{6.1. Tile cylinders and rational subspaces in~$\re^n$}
\end{center}

\bigskip

A {\em cylinder} parallel to a subspace~$J \subset \re^n$ is a
closed nonempty subset of~$\re^n$ invariant with respect to all translations parallel to~$J$. 
 In particular, if 
$J = \re^n$, then~$C = \re^n$; if $J=\{\bO\}$,  then $C$ is an arbitrary closed set. 
We usually avoid those cases and assume that~$J$ is a proper nontrivial subspace
of~$\re^n$.

We say that {\em a given tiling~$\{G+\bk\}_{\bk \in \z^n}$ contains a cylinder~$C$} if there is a subset~$E \subset \z^n$
such that~$C = \bigcup_{\bk \in E} (G+\bk)$. 
A cylinder contained in a given tiling will be referred to as~{\em a tile cylinder}. 
The set of all tile cylinders parallel to a subspace~$J$ is 
closed with respect to the operations of intersection, union,  and difference
(the latter two operations are followed by taking the closure).
Every tiling possesses  at least one cylinder parallel to~$J$, which is 
$C = \re^n$. This case is trivial. All other tile cylinders, if they exist, are called {\em proper}. 
 
A tile cylinder is {\em minimal}  if it does not contain another tile cylinder. 
 One can also specify a subspace~$J$ and consider the minimal tile cylinder
among those  parallel to~$J$. 
\begin{lemma}\label{l.75}
For an arbitrary tiling and a subspace~$J\subset \re^n$, there exists a minimal tile cylinder  parallel 
to~$J$. 
\end{lemma}  
 {\tt Proof.}  Choose a tile~$G$ and  consider the intersection~$C$ of all cylinder 
 parallel to~$J$ and containing~$G$. This cylinder may not be minimal since it can include  another cylinder~$C'$ which does not contain~$G$. In this case, however, 
 the closure of the set~$C\setminus C'$ is a cylinder containing~$G$, which contradicts  the definition  of~$C$. 
 Hence, $C$ is minimal.   
 
{\hfill $\Box$}
\medskip

\begin{lemma}\label{l.80}
All minimal tile cylinders parallel to~$J$ are integer translates  of each  other. 
Every tile cylinder parallel to~$J$ is  partitioned to the minimal ones. 
\end{lemma}  
 {\tt Proof.}  The intersection of all tile cylinders parallel to~$J$ and 
  containing the tile~$G$ is a minimal cylinder (see the proof of Lemma~\ref{l.75} above). Denote it by~$C$. Any other minimal 
  cylinder~$C'$ 
  contains some tile~$G' = G+\bk$. Hence, the cylinder~$C'- \bk$ contains~$G$ and therefore, it contains~$C$, because~$C$ is minimal. If~$C$ is  a proper subset of~$C'-\bk$, then 
 the cylinder  $C+\bk$, which contains~$G+\bk = G'$, is a proper subset of~$C'$. 
   This is impossible since~$C'$ is minimal. Thus, $C = C' - \bk$, and 
   we proved that all minimal cylinders are integer translates of each other. 
   Now we introduce an equivalence relation on the set of tiles: two tiles are equivalent if they belong to one minimal cylinder. If~$C_0$ is  an arbitrary tile cylinder 
  parallel to~$J$, then we take the maximal system of non-equivalent tiles in it. 
  Then the corresponding minimal cylinders fill~$C_0$.

{\hfill $\Box$}
\medskip 

Since $\re^n$ is also a cylinder, we obtain 
\begin{cor}\label{c.50}
For every tiling, the minimal tile cylinders form a 
partition of the space~$\re^n$. 
\end{cor}

A linear subspace~$L \subset \re^n$ of dimension~$d$ is called~{\em rational} if it is spanned by 
several integer vectors. It is known that every rational subspace has a basis of integer vectors 
that spans  the whole {\em integer part of~$L$}, which is
$L\cap \z^n$, 
by integer linear combinations. It will be referred to as an~{\em elementary basis}. 
 The characteristic property 
of the elementary basis is that its Gram matrix has the determinant one. Indeed, this condition means that the  parallelepiped spanned by the vectors of this basis 
has the $d$-dimensional volume one, and hence, it spans the integer lattice on~$L$. Every elementary basis 
can be complemented to an elementary  basis of~$\re^n$, which spans the lattice~$\z^n$. 
There is an integer matrix that maps an elementary basis to the 
first~$d$ vectors of the canonical basis of~$\re^n$. 

The linear span and the intersection of 
several rational subspaces are also  rational subspaces. This implies that 
every subspace $J \subset \re^n$ possesses the minimal by inclusion rational 
subspace~$L$ containing~$J$. This subspace will be called~{\em the minimal rational 
subspace of~$J$}.  There are countably many proper rational 
subspaces of~$\re^n$. For every vector out of those subspaces, its minimal rational subspace is the entire~$\re^n$. Thus, for almost all (in the Lebesgue measure) vectors of~$\re^n$, 
 the minimal rational subspace coincides with~$\re^n$.

\bigskip 

\bigskip 

\begin{center}
\textbf{6.2. The polynomial cones for \\ the tile B-splines}
\end{center}

\bigskip 

We consider a tile B-spline $\varphi$ of order~$\ell$, which is a convolution of 
$\ell+1$ functions~$\chi_{G}$. For this tile B-spline, $\widehat{\varphi} = \widehat{\chi_G}^{\ell + 1}$, therefore $\nill_k(\widehat{\varphi}) = \nill(\widehat{\varphi})$ for all $k \le \ell$. 
For an arbitrary~$k\le \ell$ and for every~$s=1, \ldots , q$, 
we consider the space~$P_{k,s}$
 of polynomials from~$\cP_{\Omega}$ that 
have zeros of order at least $2(k + 1)$ on the set~$\pi \{\nill (\widehat{\varphi})\, , \, J_s^{\perp}\}$ 
and the space~$P_{k} = \cP^{(k)}$  of polynomials from~$\cP_{\Omega}$ that have 
zero of order at least $2(k +1)$ at the origin. 
The following theorem claims that for a tile B-spline, 
each subspace $P_{k,s}$ coincides with~$P_k$, unless the corresponding tiling 
has a proper cylinder parallel to the minimal  rational subspace~$L_s$ of~$J_s$.   

\begin{theorem}\label{th.60}
Let~$\varphi$ be  the tile B-spline of order~$\ell$ generated by a tile~$G$;  
then for every~$k \le \ell$ and for every~$s=1, \ldots , q$ , the following holds: 
either~$P_{k, s} = P_{k}$ or the tiling~$\{G+\bk\}_{\bk \in \z^n}$ contains a proper cylinder
parallel to the minimal rational subspace of~$J_s$.   
\end{theorem}

The proof of Theorem~\ref{th.60} begins with several auxiliary results. 
For a rational subspace~$L \subset \re^n$, we say that 
its subset~$A \subset L$ is dense in~$L/\z^n$ if 
the union of the translates $A + \bk, \ \bk \in \z^n\cap L$ is dense in~$L$. 
The following fact is well-known:  
\begin{lemma}\label{l.70}
If $J \subset \re^n$ is a subspace and~$L$ is its minimal rational subspace,   
then $J$ is dense 
in~$L /\z^n$.  
\end{lemma}

\begin{prop}\label{p.80}
The minimal tile cylinder parallel to a subspace~$J$ is 
parallel to its minimal rational subspace. 
\end{prop} 
{\tt Proof.} Let $C$ be a tile cylinder and $G\subset C$ be one of its tiles, 
$C = \bigcup_{\bk \in E} (G+\bk), \, E \subset \z^n$. 
For a given subspace $V \subset \re^n$, we denote by $G_V$ the intersection 
of $G$ with the set~$\bigcup_{\bk \in \z^n} (V+\bk)$ and  
by $G_V'$ the intersection 
of $G$ with the set~$\bigcup_{\bk \in E} (V+\bk)$. 
Since the cylinder $C$ is parallel to~$J$, it follows that 
$G_J' = G_J$. Furthermore, if $L$ is the minimal rational subspace for~$J$, 
then $G_J$ is dense in~$G_L$. Hence, $G_J'$ is dense in~$G_L$. 
Consequently,  $G_L'$ is dense in~$G_L$ and 
by compactness, $G_L' = G_L$. Therefore, $C$ is parallel to~$L$.

{\hfill $\Box$}
\medskip 

Proposition~\ref{p.80} allows us to consider only rational subspaces parallel to cylinders. 
Since the linear span of several subspaces parallel to a cylinder is 
also parallel to it, there is the maximal by inclusion subspace parallel to the cylinder
 called the~{\em generatrix}. Proposition~\ref{p.80} yields that the generatrix 
 of a minimal tile  cylinder~$C$ 
 contains all (not necessarily rational) subspaces parallel to~$C$.    

\begin{cor}\label{c.70}
The generatrix of the minimal tile cylinder parallel to a subspace~$J$
is a rational subspace. 
\end{cor}

{\tt Proof of Theorem~\ref{th.60}.} As the first step, we reduce the theorem to the case~$k = 0$, i.e., 
to the tile function~$\phi = \chi_{G}$. Indeed, in view of  Proposition~\ref{p.90}, 
for the tile B-spline~$\varphi$,  either~$P_{k, s} = P_{k}$ or 
there exists a non-integer vector~$\bz \in J_s^{\perp}$ such that 
for every integer~$\bk$, the vector~$\bz+\bk$
is either orthogonal to~$J_s$ or is a zero 
of the function  $\widehat \varphi$ 
of order at least~$k+1$. 
Since 
$\widehat \varphi \, = \, \bigl(\,  {\widehat{\phi}}\, \bigr)^{\, \ell}$, it follows that 
the functions $\widehat \varphi$ and  $\widehat \phi$, have the same sets of zeros. 
So, that condition is equivalent 
to~$\widehat \phi (\bz + \bk) = 0$ for all~$\bk \in \z^n$ unless
~$\bz + \bk \in J_s^{\perp}$. Thus, one needs to prove the theorem 
for the case $k=0$, i.e., $\varphi = \phi$.

For every~$\bh \in J_s$, consider the generalized derivative~$\phi_{\bh} \in \cS'$
(the generalized derivative along the vector~$\bh$). 
We have 
$$
\widehat{\phi_{\bh}}(\bz + \bk)\ = \ 2\pi i \, (\bh, \bz + \bk)\, \widehat{\phi}(\bz + \bk)\ = \ 0 \qquad \mbox{{\rm for all}} \qquad \bk \in \z^n.  
$$
Indeed, for $\bk \ne \bO$, we have $\widehat{\phi}(\bz + \bk)\, = \, 0$, and 
for $\bk = \bO$, we have~$(\bh, \bz + \bk) \, = \, (\bh, \bz) \, = \, 0$. 
Consequently, by the Poisson summation formula, $\sum\limits_{\bq \in \z^n}
e^{2\pi i (\bz , \bq)  }\phi_{\bh}(\vardot  + \bq)\, = \, 0$. Consider the function 
\begin{equation}\label{eq.F}
F(\bx)\ = \ \sum_{\bq \in \z^n}\, 
e^{2\pi i (\bz , \bq)  }\, \phi(\bx + \bq).
\end{equation}
 The summation is well-defined because 
$\phi$ is compactly-supported. In fact, $F$ is a piecewise-constant function
corresponding to the tiling~$\{G+\bk\}_{\bk \in \z^n}$. Moreover, since the set of coefficients  
$\{e^{2\pi i ( \bz , \bq)}\}_{\bq \in \z^n}$ is bounded, it follows that 
the series converges also in~$\cS'$ and can be termwise differentiated along the direction~$\bh$. 
Therefore, $F_{\bh} = 0$ and~$F$ is an identical constant along every line 
parallel to~$\bh$. Since~$\bh$ is an arbitrary vector from~$J_s$, it follows that 
$F$ is an identical constant along every affine plane with the linear part~$J_s$. 
On the other hand, $F$ is a piecewise-constant function on the 
tiling, hence, all level sets of~$F$ are either empty or 
cylinders parallel to~$J_s$. Moreover, the vector~$\bz$ is non-integer, 
consequently, the coefficients~$e^{2\pi i ( \bz , \bq)}, \, \bq\in \z^n$, are not all 
equal. That is why the level sets of~$F$ are proper cylinders.   Now the theorem follows from Proposition~\ref{p.80}. 

Conversely, if the tiling contains a proper cylinder parallel to~$L_s$, then 
it is parallel to~$J_s$. Therefore, $\re^n$ is split to equal tile 
cylinders. Associating  an arbitrary number to each of those cylinders 
(we choose all those numbers $\{c_{\bq}\}_{\bq \in \z^n}$ to be different) we obtain a piecewise-constant function~$F(\vardot) \, = \, \sum_{\bq \in \z^n} c_{\bq} \, \phi(\cdot + \bq)$
which is constant on every cylinder. This implies that~$F$ 
is constant along any line parallel to~$J_s$, which means that 
for every~$\bh \in J_s$, the generalized derivative~$F_{\bh}$ is zero. 
Thus, $\sum_{\bq \in \z^n} c_{\bq} \, \phi_{\bh}(\cdot + \bq)\, = \, 0$. 
Therefore,  $c_{\bq} = e^{2\pi i (\bz , \bq)}$ for some 
$\bz \in \re^n$ (see~\cite{R}). Since~$F$ is not an identical constant, it follows that 
 $\bz \notin \z^n$. Now by the Poisson summation formula it follows that 
 $\, \widehat{ \phi_{\bh}} (\bz + \bk)\, = \, 0, \ \bk \in \z^n$. 
 Thus, $2\pi i \, (\bh, \bz + \bk)\, \widehat{\phi}(\bz + \bk) = 0, \ \bk \in \z^n$. 
This is equivalent to say that 
$\bz \in  \pi \, \bigl\{ \nill( \widehat \phi) \, , \,  J_s^{\perp}\bigr\}$.
Since $\nill( \widehat \varphi)= \nill( \widehat  \phi)$, we see that the set  
$\pi \, \bigl\{ \nill( \widehat \varphi) \, , \,  J_s^{\perp}\bigr\}$
contains a non-integer point, hence, $P_{k,s} \ne P_k$.

{\hfill $\Box$}
\medskip 

\bigskip 

\begin{center}
\textbf{6.3. The regularity of tile B-splines}
\end{center}

\bigskip 

Now we are formulating the main result of this section. 
\begin{theorem}\label{th.65}
For every  tile B-spline~$\varphi$, we have the  alternative:
either ${\alpha (\varphi)\, = \, \frac{1}{2}\log_{1/r} \, \rho_{\, \ell}}$, 
or the corresponding tiling~$\{G+\bk\}_{\bk \in \z^n}$ contains a proper cylinder
parallel to the minimal rational subspace~$L_1$ of~$J_1$. In the latter case, the following hold:

\smallskip 

\noindent \textbf{1)} $L_1$ is a proper subspace of~$\re^n$ invariant with respect to~$M$;  
\smallskip 

\noindent \textbf{2)} In a suitable integer 
basis of~$\re^n$, the matrix~$M$ becomes an integer upper block-triangular matrix with 
the first block~$A$ corresponding to~$L_1$: 
\begin{equation}\label{eq.blocks}
M \ = \ \left(
\begin{array}{cc}
A & *\\
0 & B
\end{array}
\right)\, , 
\end{equation}
\end{theorem} 
 {\tt Proof}. If $\alpha(\varphi) \ne \frac{1}{2} \log_{1/r} \, \rho_{\, \ell}$, then, by Theorem~\ref{th.50},  $P_1$ is a proper subspace of~$P$, which, in view of Theorem~\ref{th.60}, implies 
 that~$L_1$ is a proper subspace of~$\re^n$. 
   Since~$J_1 \subset L_1$ and $MJ_1 = J_1$, it follows that 
 the rational subspace~$ML_1$ also contains~$J_1$ and hence contains~$L_1$.
 Thus, $L_1 \subset ML_1$ and therefore, $L_1 = ML_1$ since the matrix $M$
 is nondegenerate. Choosing an  elementary basis of~$L_1$ and complementing 
it to an elementary basis in~$\re^n$, 
we obtain a basis in which~$M$  becomes an integer matrix of the form~(\ref{eq.blocks}).

{\hfill $\Box$}
\medskip

\begin{theorem}\label{th.90}
For a tile B-spline, the inequality $\alpha(\varphi) \, \ne \, \frac{1}{2} \log_{1/r}\, \rho_{\, \ell}$ is possible 
only if  
$m \ge 4$; in the  bivariate case ($n=2$), it is possible only if~$m \ge 6$. 
\end{theorem} 
{\tt Proof.} By Theorem~\ref{th.65}, the inequality 
$\alpha \ne \frac{1}{2} \log_{1/r}\, \rho_{\ell}$ implies that in a suitable elementary 
basis, the matrix~$M$ gets the upper block-triangular form~(\ref{eq.blocks})
with integer entries, 
 where the subspace~$L_1$ corresponds to  the block~$A$. 
 Since $J_1 \subset L_1$, it follows that 
 $\rho \, = \, \rho(A) > \rho(B)$. On the other hand, the matrices~$A$ and $B$
 are both  expanding and integer, hence,  their determinants are at least two
 in absolute value.   Therefore,  $m\, = \, |{\rm det}\, M|\, \ge \, 4$. In the two-dimensional case, 
 $A$ and $B$ are numbers and $|A| > |B| \ge 2$, so $|A| \ge 3$, 
 and $m = |AB| \ge 6$.

{\hfill $\Box$}
\medskip 

\begin{remark}\label{r.120}
{\em Both of the lower  bounds given by~Theorem~\ref{th.90} are sharp. 
For bivariate splines, the $6$-digit case is illustrated in Example~\ref{ex.10} 
The three-dimensional $4$-digit case is given by ``bear cylinders'' in Example~\ref{ex.20}}.  
\end{remark}

\bigskip

\section{Examples and numerical results. The super-smoothness}\label{section:numeric}

To compute the~$L_2$-regularity of a refinable 
function by Theorems~\ref{th.50h} and~\ref{th.55h} one needs to 
write the matrix of the operator~$\cT$, find  its restrictions to the 
invariant subspaces~$P_{k, s}$, compute their  spectral radii~$\rho_{k, s}$
(i.e., the Krein-Rutman leading eigenvalues) and substitute those values in formula~(\ref{eq.Hold}). 
We obtain the biggest~$k$ for which~$\varphi \in W_2^k$ and 
the H\"older exponent of the~$k$th derivative. The technical  details are given in~Appendix. 
We begin with two simple examples that prove the 
sharpness of   Theorem~\ref{th.90}. In the sequel of this section we use the same notation~$M, D$, etc.,  in different examples.  
\bigskip 

\begin{center}
\textbf{7.1.  Two examples for Theorem~\ref{th.90}}
\end{center}

\medskip

If $n\ge 3$, then the equality $\alpha(\varphi) \, = \, \frac{1}{2} \log_{1/r}\, \rho_{\, \ell}$
may fail for~$4$-digit tile B-splines; 
if $n = 2$, then it may for~$6$-digit ones. 
\begin{ex}\label{ex.10}
{\em   Consider tile B-splines~$\varphi_{\ell}$ of order~$\ell$
generated by the same diagonal dilation matrix: 
$$
M \ =  \begin{pmatrix}3 & 0 \\ 0 & 2\end{pmatrix}
$$
and by different digit systems~$D$. We have $r_1 = 3, r_2 = 2$, 
so~$r = r_1 = 3$, 
and~$m = 6$. 
For the 
system~$D \, = \, \bigl\{\, (i, j)\, : \  i\in \{0,1,2\}, \ 
j \in \{0,1\}\, \}$, the corresponding tile is a unit square, therefore, 
$\varphi_{\ell} = B_{\ell}$ and hence, $\, \alpha (\varphi_{\ell}) = \ell+\frac12$ for all~$\ell \ge 0$. 
On the other hand, the corresponding refinement equation is a direct product of two 
univariate refinement equations, hence, $\rho_{0, 1} = \frac13 , \rho_{0, 2} = \frac12$, so
$\rho_{\, 0} = \frac12$. For~$\ell =0$, we have~$\alpha = 0.5$, while~$\, 
 \frac{1}{2} \log_{1/r}\, \rho_{\, 0}\, = \, 
\frac{1}{2} \log_{1/3}\, \frac12 \, = \, 0.315...$. The same situation 
takes place for all~$\ell \ge 1$: 
$\alpha = \ell+\frac12 \, > \, \frac{1}{2} \log_{1/r}\, \rho_{\, \ell}\, = \, 
\frac{1}{2} \log_{1/3}\, 2^{-(2\ell + 1)} \, = \, \bigl(\ell + \frac12\bigr)\, \log_3 2$. 
This  anysotropic equation is a direct product of isotropic ones;  
such examples are well-known. The following example with the same matrix~$M$ is less trivial:

\begin{figure}
\begin{center}
\begin{minipage}[h]{0.35\linewidth}
\center{\includegraphics[width=1\linewidth]{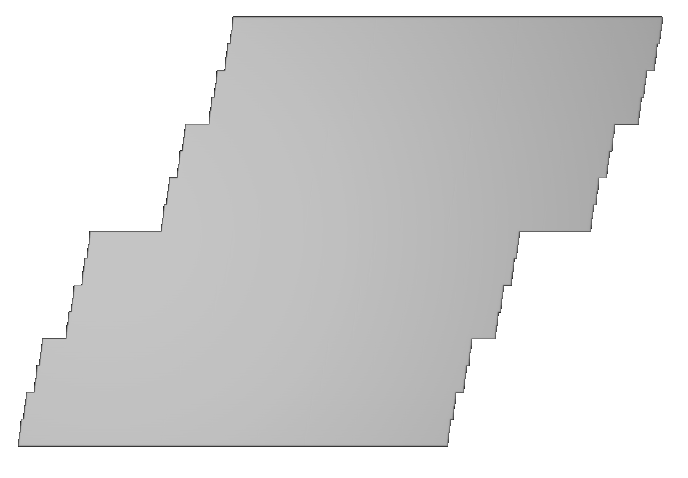}} 
\caption{{\footnotesize  6-digit tile, Example~\ref{ex.10}}}
\label{fig.ex10}
\end{minipage}
\hfill
\begin{minipage}[h]{0.35\linewidth}
\center{\includegraphics[width=1\linewidth]{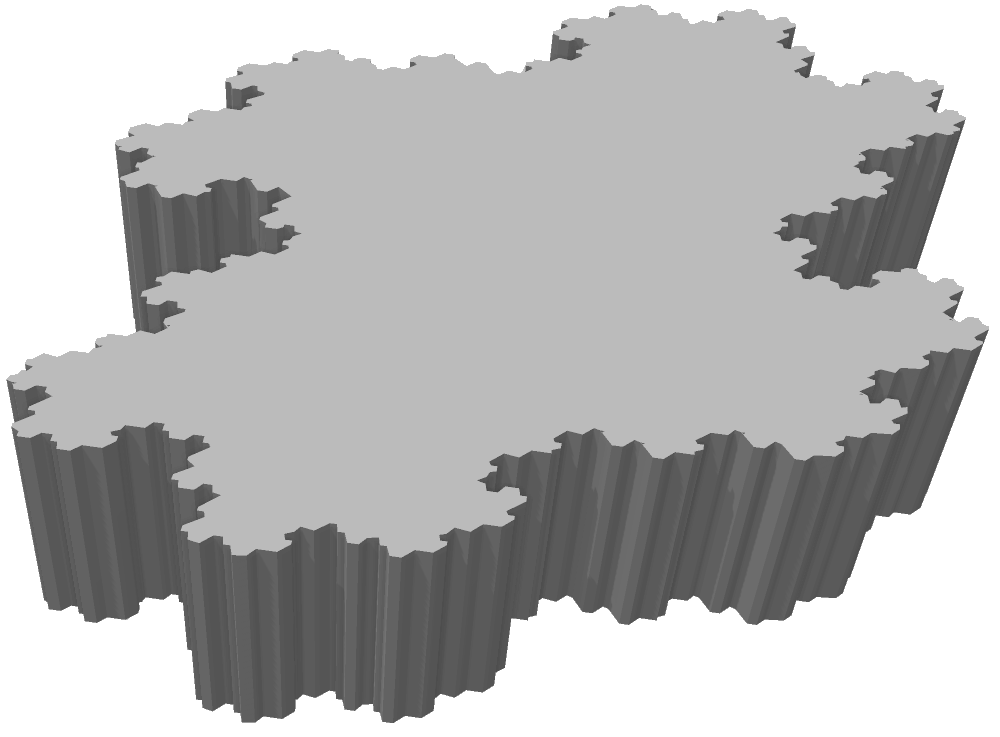}} 
\caption{{\footnotesize  4-digit bear cylinder, Example~\ref{ex.20}.}}
\label{fig.ex20}
\end{minipage}
\end{center}
\end{figure}

In the digit set~$D$,  we replace the digit~$(0,1)$ by $(3,1)$. 
The new set~$D$ is not a direct product and neither the corresponding refinement   equation. 
For~$\ell = 0$, when~$\varphi_{\ell} = \varphi$ is the indicator function 
of the tile~$G$
(fig.~\ref{fig.ex10}), we again have 
$\rho_{0,1} = \frac13, \rho_{0,2} = \frac12$ and~$\alpha = \frac12$. 
For the higher order~$\ell$, we also have~$\alpha \, > \, \frac{1}{2} \log_{1/r}\, \rho_{\ell}$ as its is seen below: 
\smallskip 

For~$\ell = 1$,   $\ \alpha =  1.630... \quad $ while $ \quad  
\frac{1}{2} \log_{1/3}\, \rho_{1}\, = \, 1.028...$; 
\smallskip 

For~$\ell = 2$, $\ \alpha =  2.261... \quad $ while  $\quad   
\frac{1}{2} \log_{1/3}\, \rho_{2}\, = \, 1.426...$; 

\smallskip 

For~$\ell = 3$,  $\ \alpha =  2.892... \quad $ while  $
 \quad  
\frac{1}{2} \log_{1/3}\, \rho_{3}\, = \, 1.824...; 
$.  
\medskip 

Theorem~\ref{th.65} guarantees the existence of the tile cylinder in this example. 
Indeed, this is the horizontal strip~$\Pi \, = \, \bigl\{ (x_1, x_2) \, : \  0\le x_2 \le 1\,  \bigr\}$ formed by  the translates~$G + (k,0)\, , k\in \z$.  
}
\end{ex} 

\begin{ex}\label{ex.20}
{\em  If~$n\ge 3$, then 
Theorem~\ref{th.90} asserts that~$\alpha \, = \, \frac{1}{2} \log_{1/r} \rho_{\, \ell}$, unless 
$m \ge 4$. The example for~$m=4$ is given  by the so-called {\em Bear cylinder} in~$\re^3$. 
The base of this cylinder is {\em Bear tile}~$G_b$ generated in~$\re^2$
by the isotropic dilation matrix: 
\begin{equation}\label{eq.bear}
M_b \ = \ \begin{pmatrix} 1 & -2 \\  1 & \quad 0\end{pmatrix}
\end{equation}
and digits~$\, D\, = \, \bigl\{\, (0,0), \, (1,0)\, \bigr\}$. 
Here~$m=2, r = \sqrt{2}$
and~$\rho_0 = 0.7607...$, therefore, $\alpha(\varphi^{(b)}) = 0.3946...$, where 
$\varphi^{(b)}$ is the indicator function of~$G_b$. The Bear cylinder~$G \, = \, [0,1]\, \times \, G_b$ satisfies a refinement equation with the 
dilation 
 $3\times 3$-matrix~$M$ that consists of  two diagonal blocks:~$\{2\}$ and $M_b$. 
 Thus, $r_1 = 2, r_2 = \sqrt{2}$, hence, $r=2$. The set of digits is~$(0,0,0), (0,1,0), (1,0,0), (1,1,0)$. 
 We have~$\alpha(\varphi) = \alpha(\varphi^{(b)}) = 0.3946...$, 
 while $\frac{1}{2} \log_{1/r} \rho_{\, 0}\, = \, 0.1973$. 
}
\end{ex} 

\medskip

\smallskip

\begin{center}
\textbf{7.2.  The super-smooth tile B-splines}
\end{center}

\smallskip

\begin{figure}[h]
\begin{center}
\includegraphics[width=0.9\linewidth]{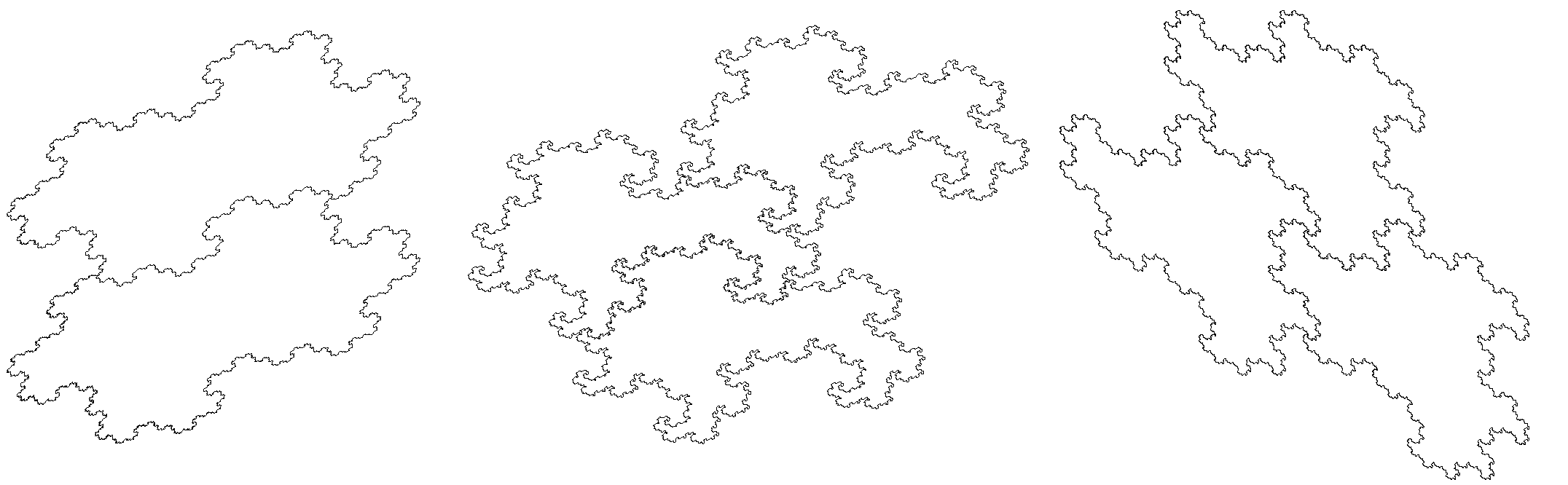}
\end{center}
\caption{Tiles generating super-smooth B-splines: the Bear  and the three-digit tiles
ThreeDig1, ThreeDig2. For each tile, the self-similar partition $G = \cup_{\bd\in D} 
M^{-1}(G + \bd)$ is shown}
\label{fig.supersmooth}
\end{figure}

We illustrate our method by computing the higher order regularity of tile B-splines 
for different dilation matrices~$M$. In particular, we find 
several families of  tile B-splines~$\varphi_{\ell}$, whose regularity exceeds that of the 
cardinal~B-splines~$B_{\ell}$, i.e.,~$\alpha (\varphi)\, > \, \ell+\frac12$.
 This phenomenon is rather surprising, we call it 
{\em super-smoothness}. Three families of tile B-splines are found in~$\re^2$:  
one of them is two-digit, i.e., $|{\rm det}\, M| = 2$
(the {\em Bear B-splines}), two others are
three-digit (see fig. \ref{fig.supersmooth}).  The results 
are shown in Tables~\ref{table0},~\ref{table1}, and~\ref{table2}. 
Table~\ref{table0} contains the H\"older exponents in~$L_2$ 
for the three ``super-smooth'' tile B-splines (lines 2 -- 4).   

\enlargethispage{\baselineskip}

The line~2 shows  the H\"older exponents of the Bear B-splines
generated by the Bear tile
 (see Example~\ref{ex.20}). This is  one of three existing 
two-digit tiles  on the plane, along with the Square and the Dragon~\cite{BG, Z2}. 
It  possesses many remarkable 
properties, see~\cite{ZP} for an overview. In particular,  their 
B-splines~$\varphi^{(b)}_{\ell}$ possesses the super-smoothness in~$C(\re)$
at least for~$\ell = 2, 3, 4$ and are useful  for construction of wavelets and 
subdivision schemes~\cite{Z1}. By the new method it is possible to 
find the regularity in~$L_2(\re)$ for higher orders~$\ell$, Table~\ref{table0}
demonstrates them for~$\ell \le 8$. We see that $\alpha(\varphi_{\ell}) \, >\,  \alpha(B_{\ell})$ for all~$\ell \ge 1$. 

The same super-smoothness phenomenon is observed for two three-digit tiles, which 
are referred to as ThreeDig1 and ThreeDig2. They are generated by the 
dilation matrices~$M_{1}, M_{2}$ given in~(\ref{eq.td}) and by the same
digit set~$D \, = \, \bigl\{ (0,0)\, , \, (1,0)\, , \, (0,1) \,\bigr\}$: 
\begin{equation}\label{eq.td}
M_{1} \ = \ 
\begin{pmatrix}1 & -2 \\ 1 & \quad 1\end{pmatrix}\ ; \qquad 
M_{2} \ = \
 \begin{pmatrix}1 & -1 \\ 1 & \quad 2\end{pmatrix}
\end{equation}
The results for the B-splines~$\varphi_{\ell}^{(2)}$ generated by ThreeDig2
are especially interesting, since their regularity are very close to the maximal 
possible level. It is well-known that if for a stable refinable function~$\varphi_{\ell}$, 
 the order of  sum rules is equal to~$\ell$, then it cannot belong to~$W^{\ell+1}_2$, i.e., 
$\alpha (\varphi_{\ell}) < \ell+1$. We see that for~$\varphi_{\ell}^{(2)}$ this 
critical level is almost achieved already for $\ell = 4$, when 
$\alpha (\varphi_{4}) > 4.9861$. For $\ell = 8$, we have  
$\alpha (\varphi_{8}) > 8.9997$. 
\smallskip 

\begin{table}[h]
\begin{center}
\begin{tabular}{c|c|c|c|c|c|c|c|c|c}
$\ell$ & 0 & 1 & 2 & 3 & 4 & 5 & 6 & 7 & 8 \\ \hline
$B_{\ell}$ & 0.5 & 1.5 & 2.5 & 3.5 & 4.5 & 5.5 & 6.5 & 7.5 & 8.5 \\
Bear &  0.3946 & 1.5372 & 2.6323 & 3.7063 & 4.7669 & 5.8173 & 6.859 & 7.893 & 8.9199 \\ 
ThreeDig1 & 0.3116 & 1.4913 & 2.6078 & 3.6918 & 4.7563 & 5.8076 & 6.8488 & 7.8819 & 8.9084 \\ 
ThreeDig2 & 0.3691 & 1.6571 & 2.8536 & 3.9518 & 4.9861 & 5.9961 & 6.9989 & 7.9997 & 8.9997 \\ 

\end{tabular}
\end{center}
\caption{The H\"older exponents of the super-smooth bivariate tile B-splines: the Bear, ThreeDig1, and ThreeDig2}
\label{table0}
\end{table}

The {\em Dragon}, which is the third two-digit plane tile, generates 
B-splines with the regularity less than the standard cardinal B-splines. 
\begin{table}[h]
\begin{center}
\begin{tabular}{c|c|c|c|c|c|c}
$\ell$ & 0 & 1 & 2 & 3 & 4 & 5   \\ \hline
Dragon & 0.2382 & 1.0962 & 1.8038 & 2.4395 & 3.0562 & 3.6688 \\
ThreeDig3 & 0.1712 & 0.9805 & 1.5872 & 2.1363 & 2.6743 & 3.21 \\
\end{tabular}
\end{center}
\caption{The H\"older exponents of bivariate tile B-splines of lower regularity}
\label{table1}
\end{table}
Among the three-digit plane tiles (there are infinitely many of them) we have found only two 
with super-smoothness. For others we have~$\alpha < \ell + \frac12$. 
Table~\ref{table1} shows the regularity of one of those examples 
(we called it ThreeDig3), which are still lower than for the Dragon tile.

We are not aware of any super-smooth tile B-splines in~$\re^3$. 
At least, there is no such examples in the two-digit case. 
In fact, 
there exist precisely seven two-digit tiles in~$\re^3$~\cite{BG}, we did the computations for all of them. 
The first one is the cube producing the cardinal B-spline~$B_{\ell}$, 
the six other produce tile  B-splines of regularity lower thar~$B_{\ell}$.  
Table~\ref{table2} shows the results for two of them (see the classification in~\cite{ZP}).  

\begin{table}[h]
\begin{center}
\begin{tabular}{c|c|c|c|c|c|c}
$\ell$ & 0 & 1 & 2 & 3 & 4 & 5   \\ \hline
Cube & 0.5 & 1.5 & 2.5 & 3.5 & 4.5 & 5.5 \\
3d Type 1 & 0.2328 & 1.3299 & 2.2125 & 2.9545 & 3.6933 & 4.4319 \\ 
3d Type 2 & 0.1173 & 0.6236 & 0.9755 & 1.3035 & 1.6296 & 1.9556 \\
\end{tabular}
\end{center}
\caption{The $L_2$-regularity  of two-digit tile B-splines in~$\re^3$}
\label{table2}
\end{table}

\begin{figure}
\begin{center}
\begin{minipage}[h]{0.35\linewidth}
\center{\includegraphics[width=1\linewidth]{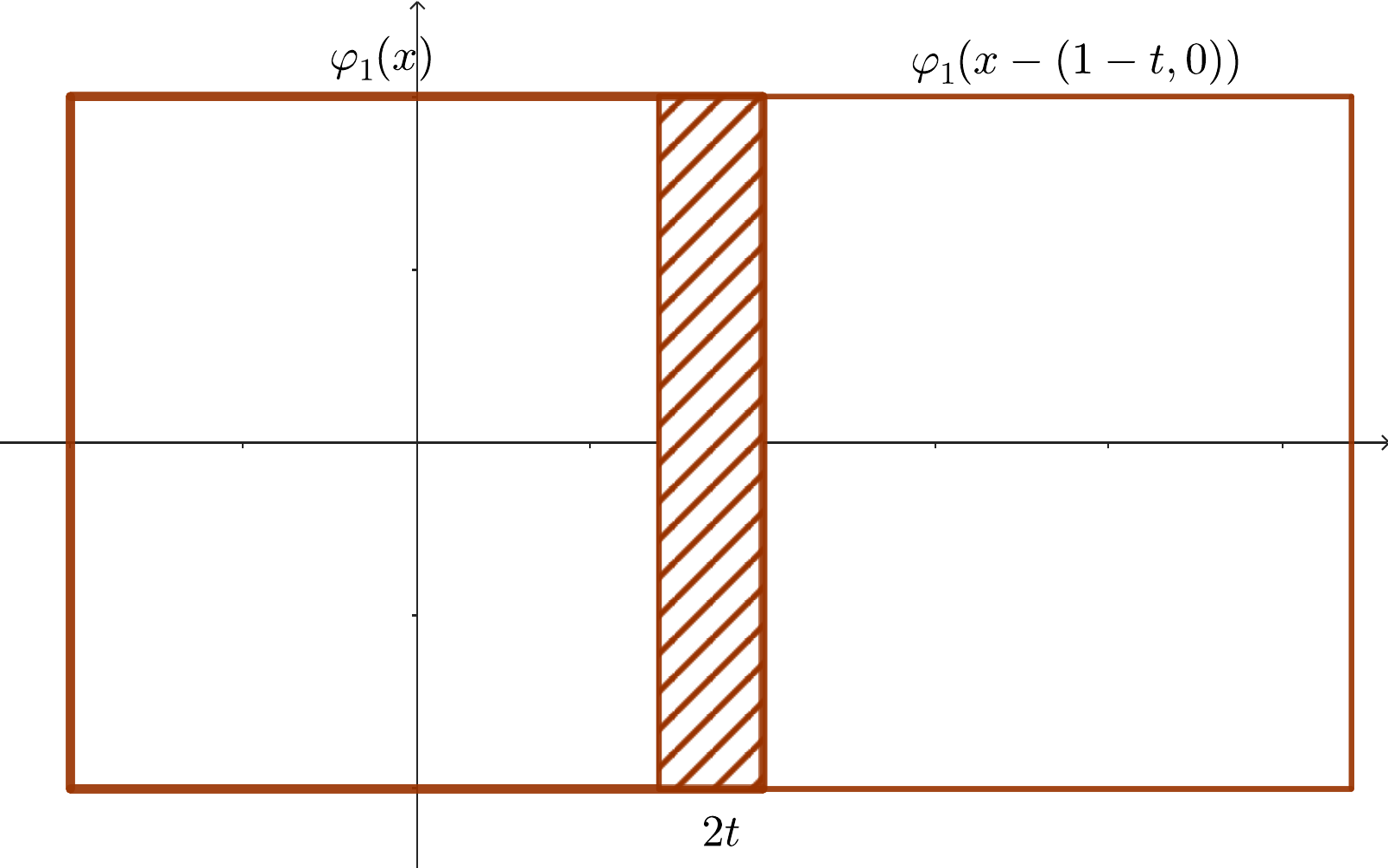}} 
\end{minipage}
\hfill
\begin{minipage}[h]{0.4\linewidth}
\center{\includegraphics[width=1\linewidth]{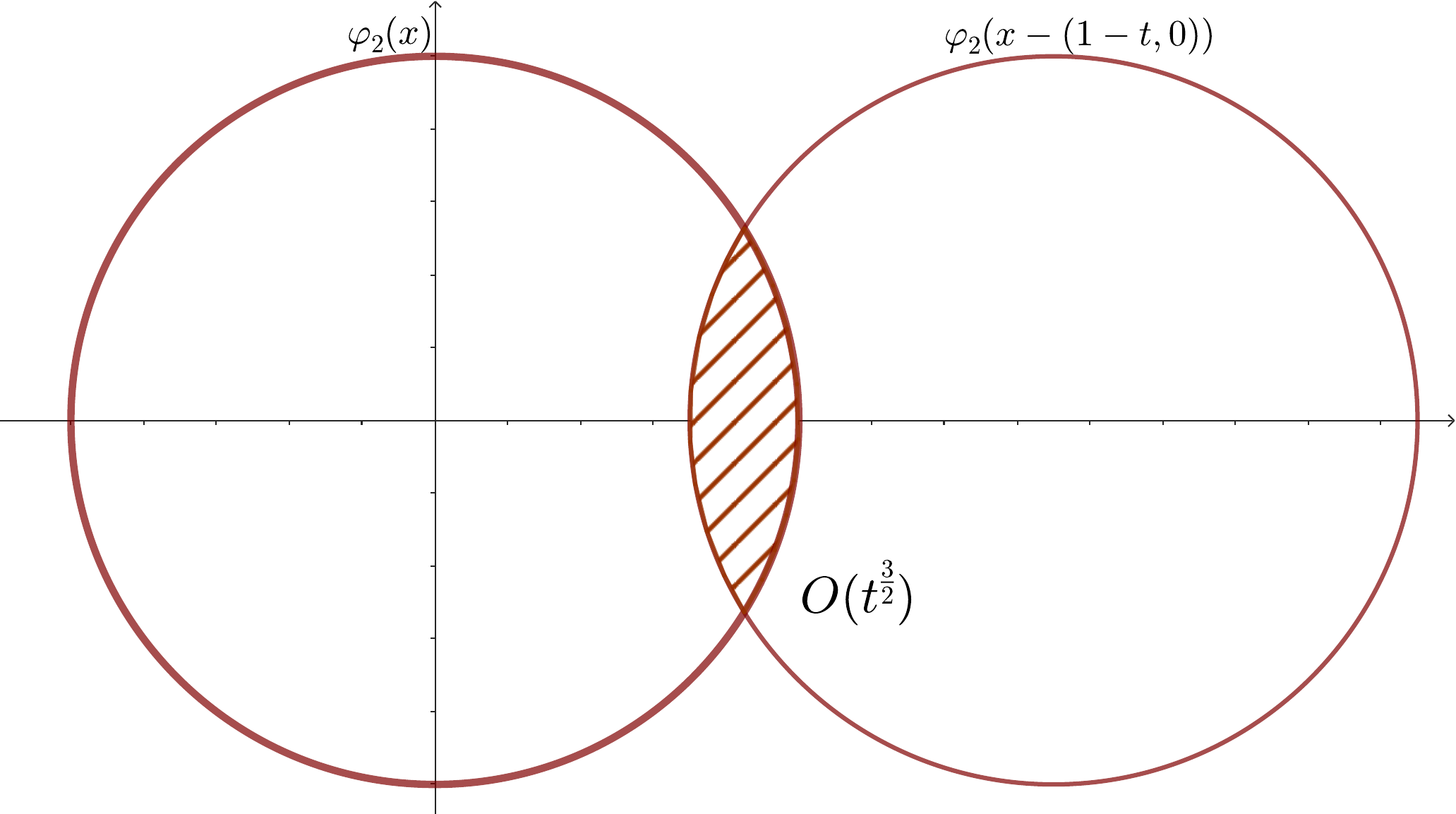}} 
\end{minipage}
\end{center}
\caption{The different smoothness of autoconvolutions of a  square and of a disc}
\label{fig.shifts}
\end{figure}

\begin{remark}\label{r.50}
{\em The super-smoothness of tile B-splines may seem contradictory: 
the autoconvolution of a square is less smooth than the autoconvolution 
(of the same order)  
of an irregular  fractal-like tile. We are not aware of any explanation of this phenomenon. 
Nevertheless, one argument can be given: 
a  convolution of a polygon with itself cannot be from~$C^1$, while 
for some other figures it can. The details are in fig. \ref{fig.shifts}. 
Let $f_1$ be the indicator of the  square~$\{(x,y)\,  : \, 
|x|\le 1, |y| \le 1 \}$ and $f_2$ be that of the unit disc centered at the origin, 
$F_i = f_i*f_i, \, i = 1,2$. Then $F_1$ is not differentiable at~$\bx = (1,0)$, while 
$F_2$ is. Indeed, $F_1(1-t, 0)$ is equal to the area of shadowed rectangle, which is 
$2t$ if $t> 0$ and $0$ otherwise. On the other  hand, it is easily shown that 
$F_2(\bx + \bh) = O(\|\bh\|^{3/2})$, hence~$F_2'(\bx) = 0$. }
   \end{remark}

\enlargethispage{\baselineskip}

\begin{center}
\textbf{7.3.  The numerical complexity}
\end{center}

\medskip

The efficiency of the Littlewood-Paley technique 
in comparison to the matrix approach  is explained by the dimension 
of the Perron eigenvalue problem. 
In the univariate case, the matrix approach deals with the matrix size
of order~$\, \frac12 \, N^2$, while the  Littlewood-Paley method does with the size 
of order~$2N$, 
where $N$ is the degree of the mask. 
 For the multivariate case, the sizes of both matrices are bigger, 
 but the size in the Littlewood-Paley method (i.e., the dimension of the 
 space~${P_{\ell}}$) is still much less than one of the matrix method.  
 Table~\ref{table3} compares these sizes for the Bear B-splines.
 The situation in other numerical examples is similar.   

\begin{table}[h]
\begin{center}
\begin{tabular}{c|c|c|c|c|c|c|c|c|c|c}
$\ell$ & 0& 1 & 2 & 3 & 4 & 5 & 6 & 7 & 8 & 9  \\ \hline
Littlewood-Paley & 4 & 11 & 20 & 55 & 76 & 101 & 162 & 203 & 242 & 333 \\
Matrix approach & 9 & 81 & 225 & 324 & 1089 & 1764 & 3969 & 5329 & 5041 & 9409 \\ 
$\alpha (\varphi_{\ell}^{(b)})$ & 0.395 & 1.537 & 2.632 & 3.706 & 4.767 & 5.817 & 6.86 & 7.893 & 8.92 & 9.94 \\
\end{tabular}
\end{center}
\caption{The matrix sizes in the two methods: computation for the  Bear B-splines.} 
\label{table3} 
\end{table}

\bigskip

\section{Optimality  of the tile B-splines}\label{section:optimality}

We establish the following optimal property of tile B-splines: 
they possess the smallest number of nonzero coefficients 
among all refinable functions of a given order of approximation. 

For the sake of simplicity, we  restrict to the case~$m=2$, i.e., 
to two-digit refinable functions. Thus, 
among all refinable functions corresponding to a given 
dilation matrix~$M$ with $|{\rm det}\, M| = 2$  and 
generating (by integer shifts) the space of algebraic polynomials of a given degree,  
the tile B-spline has the minimal number~$|Q|$ of nonzero coefficients. 
This property, which  will be referred to as 
{\em the minimal support of the mask}, also holds for  tile 
B-splines  in the following classes of refinable functions: 
\smallskip 

1) of a given 
rate of approximation of smooth functions by contractions of integer translates; 
\smallskip 

This implies that among all refinable functions generating the multiresolution analysis
(MRA)
$\{V_j\}_{j\in \z}$  with a given 
 rate of approximation by the subspaces~$V_j$, the B-splines have the minimal 
 support of the mask.  Consequently, they have the fastest algorithm of computation 
 of the wavelet coefficients produced by the corresponding MRA. 
\smallskip 

2) of a given degree of convergence of the subdivision scheme. 
\smallskip

This  means, in particular,  that the subdivision scheme 
generated by the tile B-spline has the least  algorithmic complexity 
(the number of operations to realize one iteration) among all 
schemes reproducing the space of polynomials of a given degree.  
\medskip

The digit set for the 
matrix~$M$ is 
$D \, = \, \{\bO, \bd\}$, where $\bd \in \z^n \setminus M\z^n$.
The digit set for the 
matrix~$M^T$ is 
$D^* \, = \, \{\bO, \bd^*\}$, where $\bd^* \in \z^n \setminus M^T\z^n$.

\smallskip 

The following theorem is well-known (see, for instance,~\cite{KPS}). 
 \medskip

\noindent \textbf{Theorem~B}. {\em For a given integer~$\ell \ge 0$,  the following properties of a stable two-digit refinable function 
are equivalent: 
\smallskip 

1) the integer translates $\{\varphi (\vardot + \bk)\}_{\bk \in \z^n}$
generate the space of algebraic polynomials of degree~$\ell$; 
\smallskip 

2) the refinement equation satisfies the sum rules of order~$\ell$, i.e., 
the mask~$\bc(\bxi)$ has zero of order $\, \ge \, \ell + 1$ at the point~$\bxi = M^{-T}\bd^*$. 
\smallskip 

3) the function~$\varphi$ has the rate of approximation at least $\ell+1$, 
i.e., the distance from an arbitrary  function~$f \in C^{\ell+1}(\re^n)$
to the space generated by the functions~$\{\varphi (N \vardot + \bk)\}_{\bk \in \z^n}$
is $O(N^{-\ell - 1})$. }
\medskip

We denote by $\nu(\varphi)$ the {\em accuracy} of the refinable function~$\varphi$, i.e., 
the maximal $\ell \ge 0$ for which the conditions 
1) - 3) from  Theorem~B are satisfied. This is known that if $\varphi \in W^{\ell}_2$, then 
$\nu(\varphi) \ge \ell$~\cite{KPS}.  
\medskip 

\begin{theorem}\label{th.100}
If a two-digit refinable function~$\varphi$ satisfies~$\nu (\varphi) \ge \ell$, 
then its refinement equation possesses at least~$\ell+2$ nonzero coefficients. 
\end{theorem}  
The mask of the tile B-spline~$\bc(\bxi)\, = \, 
\bigl(\frac{e^{2\pi i \, (\bd, \bxi)}+ 1}{2} \bigr)^{\ell+1}$ has 
exactly $\ell+2$ nonzero coefficients. Therefore, for each dilation matrix~$M$ such that~$|{\rm det}\, M|= 2$, among all refinable functions with~$\nu(\varphi) \ge \ell$, the tile B-spline 
of order~$\ell$ has the minimal support of the mask~$|Q|=\ell + 2$. 
\bigskip 

\noindent {\tt Proof of Theorem~\ref{th.100}}. Consider an arbitrary 
refinement equation with a mask~$\bc(\vardot)$ satisfying the sum rules of order~$\ell$. 
 Denote $\bv = M^{-T}\bd^*$ 
and define  the univariate polynomial~$\bq(t) = \bc\, (2t\, \bv), \ t\in \re$. 
This polynomial
has a root of multiplicity at least~$\ell + 1$ at the point~$t=\frac12$. 
Since $|{\rm det}\, M|=2$, it follows that $2 \, M^{-T}$ is an integer matrix, and therefore, 
$2\bv \in \z^n$. This implies that the scalar products $(2\bv, \bk)$
are integer for all~$\bk \in \z^n$ and hence, 
 $$
 \bq (t) \ =\ \frac12\, \sum_{\bk\in \z} \, c_{\bk}\, e^{2\pi i \, (2t\bv , \bk)}\ = \ 
  \sum_{s \in \z} 
  \left( \, \sum_{\bk \in \z^n, \, (2\bv , \bk) = s} \, \frac12\, c_{\bk}\, \right)\, 
  e^{2\pi i \, st} \ = \ \sum_{s \in \z} \, q_s \, e^{2\pi i \, st}, \
 $$
where $\ q_s \, = \, \frac12\, \sum\limits_{\bk \in \z^n, (2\bv , \bk) = s} \,  c_{\bk}$. 
Since $\bq^{(r)}\bigl(\frac12 \bigr)\, = \, 0$ for all~$r=0, \ldots , \ell$, 
we have 
$$
 \sum_{s \in \z} \, q_s \bigl(2\pi i s \bigr)^r e^{\pi i s} \quad = \quad 
0, \qquad r=0, \ldots , \ell\, . 
$$
In other words, 
\begin{equation}\label{eq.moments}
 \sum_{s \in \z} \, q_s \, s^r\, (-1)^s  \quad = \quad 
0, \qquad r=0, \ldots , \ell\, . 
\end{equation}
Observe that $\sum_{s} q_s \, = \, \bq(\bO)\, = \, 1$. For an arbitrary $h \in \re$, we 
consider the vector of moments~$\tilde \bh = \bigl(1, h, \ldots , h^{\ell} \bigr) \in \re^{\ell + 1}$. Every $\, N \le \ell +1$ different  vector of moments are linearly independent, since they can be complemented  to $\ell+1$ columns of a Vandermonde matrix. Assume that among the numbers 
$\bigl\{ s\, =\, (2\bv , \bk) \, : \ \bk \in \z^n\, , c_{\bk} \ne 0 \bigr\}$ there are 
exactly  $N$ different ones. Denote them by~$s_1, \ldots , s_N$. 
Then equation~(\ref{eq.moments}) reads $\sum_{j =1}^N (-1)^{s_j}\, q_{s_j} \, \tilde \bs_j \, = \, 0$. Therefore, the number~$N$ of vectors~$\tilde \bs_j$ is at least~$\ell +2$. 
On the other hand, each nonzero coefficient~$c_{\bk}$ generates exactly one 
number~$s = (2\bv , \bk)$ (which may coincide for different~$\bk$), hence, there are at least $N \ge \ell +2$ nonzero 
coefficients~$c_{\bk}$.

{\hfill $\Box$}
\medskip 

\begin{remark}{\em 
The tile B-spline is not the only optimal refinable function. One may classify all  refinement 
 equations with~$\nu (\varphi) \ge \ell$ and with $\ell +2$ nonzero coefficients. 
   We take arbitrary~$\ell + 2$ numbers~$s_1, \ldots , s_{\ell+2}$, not of the same 
   parity, and for each~$i=1, \ldots , \ell+2$,  find 
   an integer point~$\bk_i$ such that~$(2\bv , \bk_i) = s_i$. 
   This is possible because $\bv \notin \z^n$, and hence~$2\bv \notin 2\z^n$. 
   Then find the solution 
   $q_{s_1}, \ldots , q_{s_{\ell+2}}$ of the equation~(\ref{eq.moments}).
   We choose  the parities of $s_i$ so that $\sum_i q_i \ne 0$. 
   Normalizing the solution as $\sum_i q_i  = 1$ and taking~$c_{\bk_i} = 2q_i$
   we obtain the desirable refinement equation.

}
\end{remark}

\bigskip

\bigskip

\section{Tile  subdivision schemes}\label{section:subdivisions}

We restrict ourselves to subdivision schemes on rectangular grid. 
From the results of Section~\ref{section:optimality} we obtain
the following optimal property of schemes generated by tile B-splines.  
\begin{theorem}\label{th.110}
If a two-digit subdivision scheme preserves the space of algebraic polynomials 
of degree~$\ell \ge 0$, then it has at least~$\ell + 2$ coefficients. 
\end{theorem}
Thus, the tile subdivision schemes possess the minimal possible 
number of coefficients among all schemes with a given rate of approximation. 
Moreover, some of them are super-smooth. The following result was obtained in~\cite{Z1}.  
\smallskip 

\noindent \textbf{Theorem~C.}
{\em Let~$\varphi_{\ell}^{(b)}$ be the B-spline generated  by the Bear tile. Then 
$\varphi_{2}^{(b)} \in C^2(\re^2)$ and $\varphi_{3}^{(b)} \in C^3(\re^2)$}. 
\smallskip 

Note that the cardinal B-splines of orders~$2$ and $3$ do not belong to
$C^2$ and to~$C^3$ respectively. In view of Theorem~\ref{th.110}, this means that
 $\varphi_{2}^{(b)}$ has the minimal number (in this case, four) of coefficients among all~$C^2$ refinable 
 functions and the same is for $\varphi_{3}^{(b)}$ among~$C^3$ functions
 (five coefficients). 
The corresponding tile subdivision schemes involve the minimal number of points 
of the grid ($4$ and $5$ respectively) in each iteration, among all~$C^2$ and 
$C^3$ schemes.

\begin{figure}
\begin{center}
\begin{minipage}[h]{0.35\linewidth}
\center{\includegraphics[width=1\linewidth]{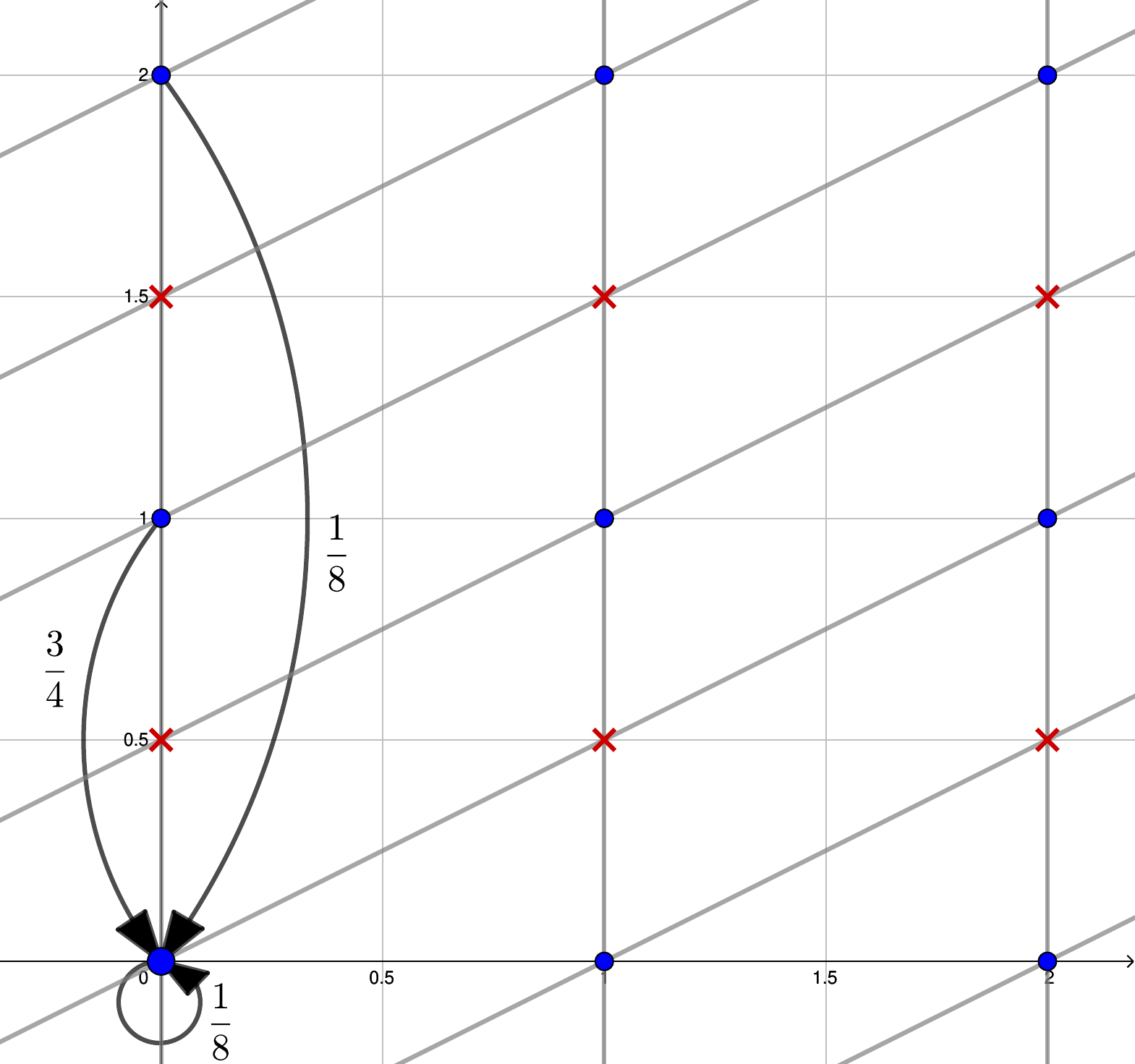}} 
\end{minipage}
\hfill
\begin{minipage}[h]{0.35\linewidth}
\center{\includegraphics[width=1\linewidth]{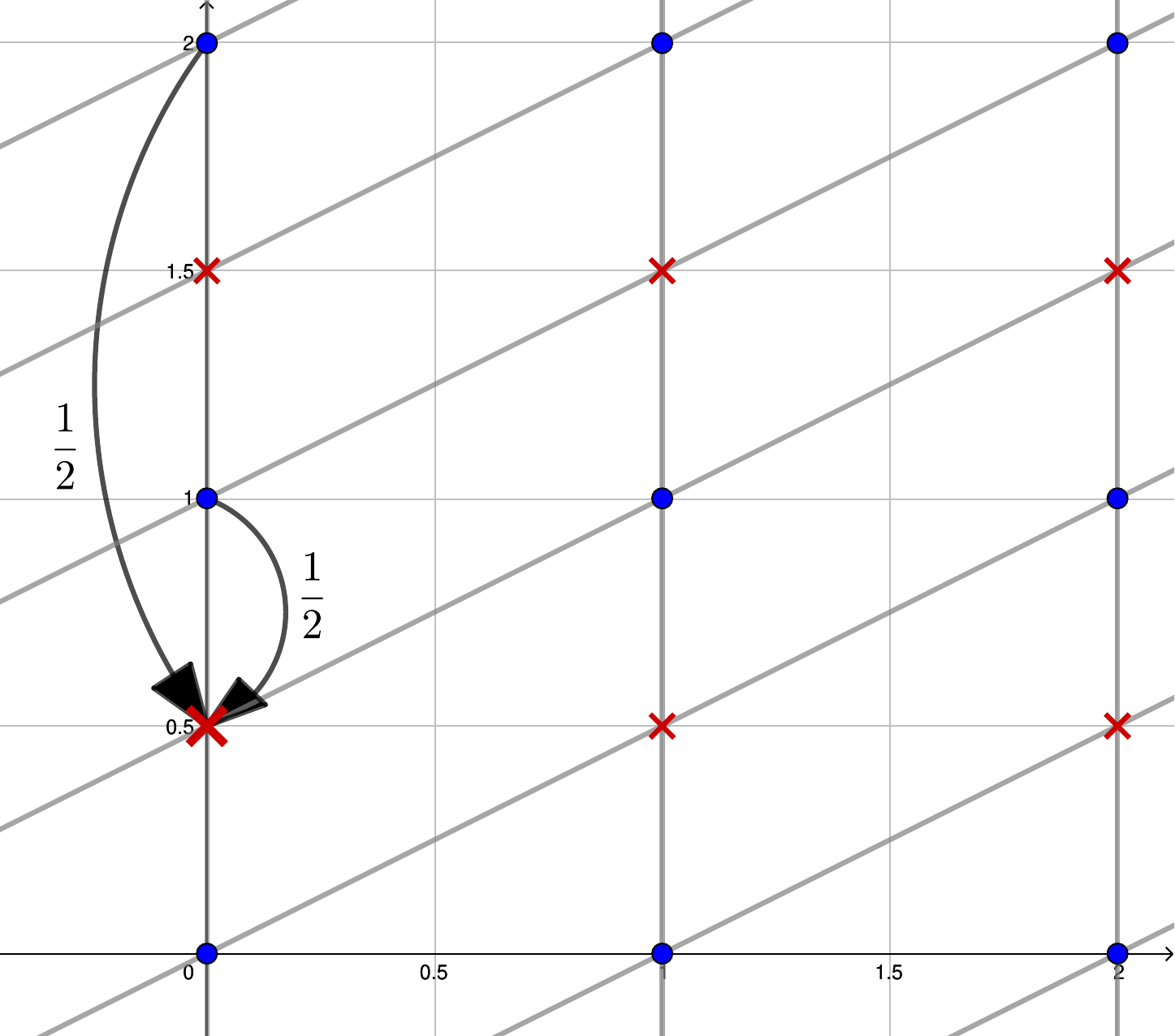}} 
\end{minipage}
\end{center}
\caption{First iteration of the $\varphi_{3}^{(b)}$ scheme}
\label{fig.scheme}
\end{figure}

\begin{figure}
\begin{center}
\begin{minipage}[h]{0.35\linewidth}
\center{\includegraphics[width=1\linewidth]{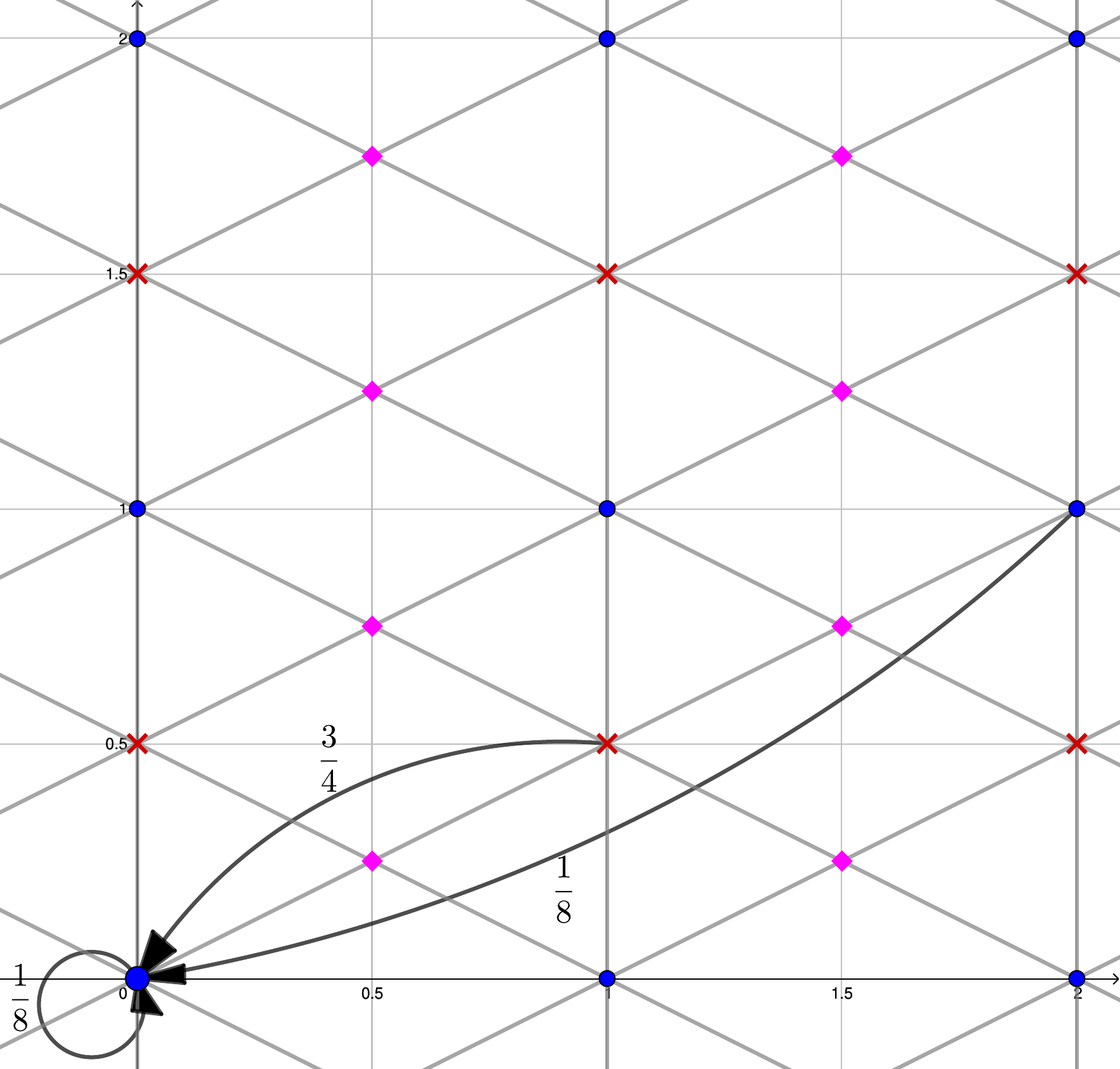}} 
\end{minipage}
\hfill
\begin{minipage}[h]{0.35\linewidth}
\center{\includegraphics[width=1\linewidth]{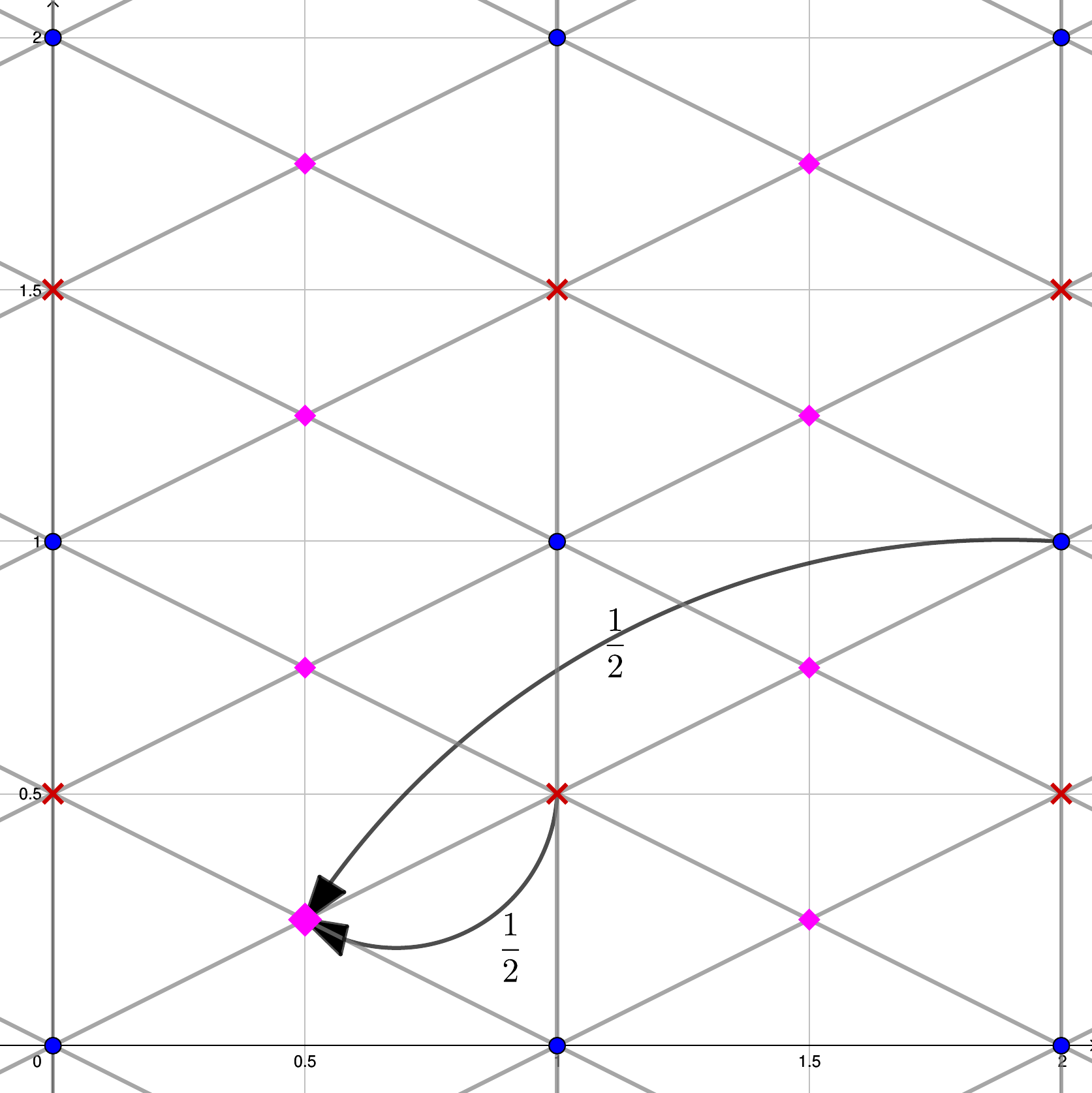}} 
\end{minipage}
\end{center}
\caption{Second iteration}
\label{fig.scheme1}
\end{figure}

The coefficients of the mask of the Bear scheme are located on one line: 
$c_{\bk},\, \, \bk = (s, 0), \, s = 0, \ldots , \ell+1$. 
The dilation matrix~$M_b$ given by~(\ref{eq.bear}) 
defines a composition of the $\sqrt{2}$-expansion and the elliptic rotation by the 
angle~$\beta = {\rm atan}\, \sqrt{7}$. Since~$\beta$ is irrational modulo~$\pi$, 
the images of the integer grid~$M^{-k}\z^2$ produced with iterations of the subdivision scheme 
covers an everywhere dense set of directions. The first iterations of 
the~$\varphi_{3}^{(b)}$ subdivision scheme ($C^3$ scheme with five coefficients)
are shown in figures~\ref{fig.scheme},~\ref{fig.scheme1}.

\bigskip 


\bigskip 

\bigskip








\bigskip 

\begin{center}
\large{\textbf{Appendix}}
\end{center}

\bigskip

{\tt Proof of Lemma~\ref{l.50}.} Clearly,  the operator~$\cT$ respects the positivity of the function. 
Therefore, if the space of polynomials~$P_{\ell, s}$ is invariant with respect to~$\cT$, then so is the cone $K_{\ell, s}$. It is well-known that the operator~$\cT$ respects 
the space of polynomials~$\cP_{\Omega}$. 
The space~$P_{\ell, s}$ is also invariant,   provided the sum rules of order~$\ell$ are satisfied. Indeed, if a polynomial~$\bp$ has zero of order 
at least $2(\ell+1)$ at the point $\bxi = \bO$, then so does the polynomial~$\cT\, \bp\, $: 
for all digits~$\bd^* \ne \bO$, the 
function~$\bigl|\bc \, \bigl(M^{-T}(\bxi + \bd^* )  \bigr)  \bigr|^2$ has zero of 
order at least~$2(\ell+1)$ at the point~$\bxi = \bO$ (because of the 
sum rules), and for~$\bd^* =  \bO$ so does the function 
$\bp \, \bigl(M^{-T}(\bxi + \bd^* )\bigr)$ (because of the 
assumption). Hence, all terms in the sum~(\ref{eq.T0})
have zero of order at least 
$2(\ell+1)$ at the point $\bxi = \bO$. 

It remains to show that if $\bp$ has a root~$\bz \in  \pi \, \bigl\{ \nill( \widehat\varphi) \, , \,  J_s^{\perp}\bigr\}$ of order at least~$2(\ell + 1)$, 
then so has the polynomial $\cT\bp$. 
To this end we consider arbitrary vectors~$\bh_i \in J_s, \, i = 0, \ldots , \ell$,  and 
denote by $g = \varphi_{\bar \bh}\, = \, 
\frac{\partial}{\partial \bh_0 \cdots  \partial \bh_{\ell}} \, \varphi$ the mixed 
derivative of order~$\ell + 1$ of the refinable function~$\varphi$ along
the directions~$\bh_i$.  
The function  $\Phi_g$ is a  polynomial from~$\cP_{\Omega}$ which 
has a root of order at least~$2(\ell + 1)$ at the point~$\bz$.  
Indeed, $\widehat{g}(\bxi) = \widehat{ \varphi_{\bh}}(\bxi)\, = \, 
2\pi i (\bh , \bxi) \widehat{ \varphi }(\bxi)$. 
If an integer shift~$\bz + \bk$ belongs to~$\nill_{\ell}( \widehat\varphi)$, then 
this point is a root of $ \widehat{ \varphi }$ of order at least~$\ell+1$; 
 if $\bz + \bk \in J_s^{\perp}$, 
then $(\bh_i , \bz + \bk)= 0$ for all~$i$, and we again 
have a zero of order at least~$\ell+1$ of the function~$\widehat{g}$ 
at the point~$\bz + \bk$. 
Hence,~$\Phi_g$ have zeros of order at least~$\ell+1$ at all the points~$\bz + \bk, \ \bk \in \z^n$.
This is true for every~$z \, \in \, 
\nill_{\ell}( \widehat\varphi) \, \cup \,  J_s^{\perp}$, therefore, 
 $\Phi_g$ vanishes at the set~$\pi \, \bigl\{ \nill_{\ell}( \widehat\varphi) \, , \,  J_s^{\perp}\bigr\}$.  On the other hand,~$\cT (\Phi_g) \, = \, \Phi_{\cF g}$, 
where the function $\cF g \, = \, \cF (\varphi_{\bar \bh})\, = \, \varphi_{M\bar \bh}$. 
Therefore, $\cT (\Phi_g) \, = \, \Phi_{\varphi_{M\bar \bh}}$  with $M\bar \bh \in J_s^{\, \ell}$, because~$J_s$ is 
an invariant subspace of~$M$. Hence, the polynomial~$\cT (\Phi_g)$
also vanishes on the 
set~$\pi \, \bigl\{ \nill_{\ell}( \widehat\varphi) \, , \,  J_s^{\perp}\bigr\}$. 
It remains to note that  zeros of the polynomial~$\cT \bp$ depend only 
on the location and orders of zeros of~$\bp$. Hence, if a polynomial vanishes on~$\pi \, \bigl\{ \nill_{\ell}( \widehat\varphi) \, , \,  J_s^{\perp}\bigr\}$, then so does~$\cT \bp$, which completes the proof.

{\hfill $\Box$}
\medskip 

\begin{center}
\textbf{Algorithms and implementation details}
\end{center}

To write the matrix of~$\cT$ we need 
 to find the set $\Omega = \bigl(Y - Y\bigr) \cap \z^n$, which can be obtained by 
 a special iterative algorithms, similar to ones 
 presented in~\cite{CM, CP}.  We give a brief description below.

Using the definition of $Y$ (subsection~5.1), we obtain  
$$
Y - Y\quad =\quad \left\{ \  \sum_{j =1}^{\infty} M^{-j}\, (\bs_j - \bt_j) \ : \quad \bs_j, \bt_j \in Q\, , \ j \in \n \ \right\}\, .
$$ 
Denote  $Q' = Q - Q$ and rewrite:  
$$
Y - Y\quad =\quad \left\{ \  \sum_{j =1}^{\infty} M^{-j}\, \bs_j \ : \quad \bs_j \in Q'\, , \ j \in \n \ \right\}\, .
$$ 
Define the map~$\cJ$ on the set of compact subsets~$X \subset \re^n$ as follows: 
$\cJ\, X \, = \, M^{-1}(X + Q')$. A set $X$ is {\em invariant with respect to $\cJ$}, if $\cJ \, X \, \subset \, X$. Similarly we define a map~$\cJ_0$  on the set of finite 
subsets of integers~$X \subset \z^n$: 
$\cJ_0\, X \, = \, \cJ X \cap \z^n$. The invariance  with respect to $\cJ_0$ is defined in the same way.

\begin{prop}\label{p.65} 
The set $\Omega$ is invariant with respect to $J_0$. 
\end{prop}
{\tt Proof.}
Let $\bz$ be an arbitrary element of the set $\cJ_0 \, \Omega$, i.e., 
$\bz = M^{-1}(\bs + \bq')$, $\bz \in \z^n$ for some $\bs \in \Omega, \bq' \in Q'$. 
Since $\bs \in \Omega$, we have $\bs = \sum_{j =1}^{\infty} M^{-j}\, \bs_j$, $\bs_j \in Q'$. Therefore $\bz = M^{-1}(\bs' + \bq') = \sum_{j =1}^{\infty} M^{-j}\, \bs_j'$, where $\bs_1' = \bq'$, $\bs_j' = \bs_{j - 1}$ for $j \ge 2$. Thus, $\bz \in Y - Y$. Since $\bz \in \z^n$, it follows that $\bz \in \Omega$, which completes the proof. 

{\hfill $\Box$}
\medskip 

To find the set~$\Omega$, we start with a subset  $\Omega_0 \subset \z^n$, 
which will be specified later,  
then  define $\Omega_{i + 1} \, = \, \Omega_i \, \cup \, \cJ_0 \Omega_i$ 
for $i \ge 0$. Thus, we obtain a non-increasing (by inclusion) sequence~$\{\Omega_i\}_{i\ge 0}$. 
Since $M$ is an expanding matrix, it follows 
that there exist~$\varepsilon > 0$ and an ellipsoid~$E\subset \re^n$ (Lyapunov ellipsoid) such that~$M^{-1}E 
\subset (1-\varepsilon)E$. Hence, for sufficiently 
large~$\lambda$, we have~$M^{-1}(\lambda E + Q') \, \subset \lambda E$. 
This implies that the sequence~$\{\Omega_i\}_{i\ge 0}$ is bounded above: 
all sets~$\Omega_i$ are inside $\lambda E$, where $\lambda$ is large enough. 
Therefore, for some $k$, we have $\Omega_{k+1} = \Omega_k$. 
It remains to put~$\Omega = \Omega_k$. 

The set $\Omega_0$ can be defined as
$H+\Delta$, where $\Delta$ is the set of integer
vectors~$\bx$ such that $\sum_{i=1}^n |x_i| \le 1$, 
and $H$ is an arbitrary compact set such that~$\cJ H \, \subset \, H$~\cite{P20}. 
For example, one may take the aforementioned Lyapunov ellipsoid~$E$, 
which can be found by solving the corresponding  semidefinite programming problem.

\bigskip

\textbf{Acknowledgements.} The authors are grateful to S.Tikhonov and I.Shkredov for useful discussions.


\end{document}